\documentclass
{amsart}%
\usepackage{hyperref}
\usepackage{cite}
\usepackage{amssymb}
\usepackage{amsmath}
\usepackage{amsthm}
\usepackage{amsfonts}
\usepackage{graphicx}%
\usepackage{comment}

\usepackage{color}
\usepackage{bbm}

\usepackage[toc,page]{appendix}
\usepackage{verbatim}

\newtheorem{theo}{Theorem}[section]
\newtheorem{pro}[theo]{Proposition}
\newtheorem{lem}[theo]{Lemma}
\newtheorem{cor}[theo]{Corollary}
\newtheorem{con}[theo]{Conjecture}
\newtheorem{defin}[theo]{Definition}

\theoremstyle{definition}
\newtheorem{exa}[theo]{Example}

\newtheorem{notation}[theo]{Notation}

\theoremstyle{remark}
\newtheorem{rem}[theo]{Remark}

\numberwithin{equation}{section}
\newtheorem*{acknowledgement}{Acknowledgement}

\newcommand{\SU}{\operatorname{SU}}
\newcommand{\su}{\mathfrak{su}}
\newcommand{\Ad}{\operatorname{Ad}}
\newcommand{\tr}{\operatorname{tr}}

\newcommand{\Ric}{\mathrm{Ric}}

\begin{document}

\newcommand{\sgn}{\operatorname{sgn}}

\newcommand{\pminterval}[1]{\Biggl[- {#1}, {#1}\Biggr]}

\title[$\SU(2)$]{Left-invariant geometries on $\SU(2)$ are uniformly doubling}

\author[Eldredge ]{Nathaniel Eldredge {$^{\dag}$}}
\thanks{\footnotemark {$^{\dag }$} Research was supported by a grant from the Simons Foundation (\#355659, Nathaniel Eldredge).}
\address{$^{\dag }$ School of Mathematical Sciences\\
  University of Northern Colorado\\
  Greeley, CO 80639, U.S.A.
}
\email{neldredge@unco.edu}

\author[Gordina]{Maria Gordina{$^{\dag \dag}$}}
\thanks{\footnotemark {$\dag \dag$} Research was supported in part by NSF grant DMS-1007496 and the Simons Fellowship.}
\address{$^{\dag \dag}$ Department of Mathematics\\
University of Connecticut\\
Storrs, CT 06269,  U.S.A.}
\email{maria.gordina@uconn.edu}

\author[Saloff-Coste]{Laurent Saloff-Coste{$^{\ddag }$}}
\thanks{\footnotemark {$\ddag $} Research was supported in part by NSF grant DMS-1404435 and NSF grant DMS-1707589.}
\address{$^{\ddag }$
Department of Mathematics \\
Cornell University\\
Ithaca, NY 14853, U.S.A
}
\email{lsc@math.cornell.edu}





\begin{abstract} A classical aspect of Riemannian geometry is the study of estimates that hold uniformly over some class of   metrics. The best known examples are eigenvalue  bounds under
 curvature assumptions. In this paper, we study the family of all left-invariant geometries on $\SU(2)$. We show that left-invariant geometries on $\SU(2)$ are uniformly doubling and give a detailed estimate of the volume of balls that is valid for any of these geometries and any radius.  We discuss a number of consequences concerning the spectrum of the associated Laplacians and the corresponding heat kernels.
\end{abstract}

\maketitle \tableofcontents

\renewcommand{\contentsname}{Table of Contents}

\section{Introduction}\label{s.intro}

\subsection{A conjecture and the main result}\label{s.MainResult}
This work is devoted to the uniform analysis of the family of all left-invariant Riemannian metrics on the Lie group $\SU(2)$. This is the simplest case of a natural problem we now describe.

Let $K$ be a connected real compact Lie group, and let $\mathfrak{L}(K)$ denote the family of all left-invariant Riemannian metrics $g$ on $K$. We conjecture that for each group $K$, many aspects of spectral analysis of the corresponding Laplace-Beltrami operator $\Delta_{g}$ with $g \in \mathfrak{L}(K)$, as well as the analysis of the associated heat equation $\partial _t-\Delta_{g}=0$, can be controlled \emph{uniformly} over $\mathfrak{L}(K)$. Recall that the operator $-\Delta_{g}$ has non-negative discrete spectrum with finite multiplicity, and so we can consider  the lowest non-zero eigenvalue denoted by $\lambda_{g}$.

It was shown in \cite{LiP1980a} (see also \cite{JudgeLyons2017}) that
on any compact homogeneous manifold, one has the lower bound
\begin{equation}\label{li-eigenvalue-lower}
\lambda_g \geqslant \frac{\pi^2}{4 \operatorname{diam}_{g}^{2}}.
\end{equation}
We conjecture that a matching upper bound holds uniformly over $g \in
\mathfrak{L}(K)$, so that
\begin{equation}\label{e.1.0}
  \lambda_{g} \leqslant \frac{C_K}{\operatorname{diam}_{g}^{2}}
\end{equation}
where the constant $C_K$ may depend on $K$ but not on $g$.

In terms of the heat equation, we conjecture that there are constants
$c_i=c_i(K)\in (0,\infty)$, $i=1, \dots, 4$ such that the fundamental
solution (heat kernel) $(t,x,y)\mapsto p_{t}^{g}(x, y)$ of the heat
equation on $(K,g)$ satisfies
\begin{equation}\label{e.1.1}
\frac{c_1}{V_{g}(\sqrt{t})}\exp(-c_2 d_{g}(x, y)^2/t)
\leqslant  p_{t}^{g}(x, y) \leqslant \frac{c_3}{V_{g}(\sqrt{t})}\exp(-c_4 d_{g}(x,y)^2/t).
\end{equation}
Here $V_{g}(r)$ denotes the volume of the ball of radius $r$ with
respect to the Riemannian volume measure $\mu_{g}$; $d_{g}(x,y)$
denotes the Riemannian distance between $x$ and $y$; and
$\operatorname{diam}_{g}$ denotes the diameter of $K$ with respect to
$d_g$.

One reason to believe that this conjecture might be true is that it
can be reduced to a simpler question. Let $\left( X, d, \mu
\right)$ be a metric measure space, that is, $\left( X, d \right)$ is
a metric space and $\mu$ is a Borel measure on $X$. By $B_{r}\left(x
\right)$ we denote the ball centered at $x \in X$ of radius $r > 0$
with respect to the distance $d$. The metric measure space  $\left( X, d, \mu \right)$ is \textbf{volume doubling} if
\begin{equation}\label{e.1.2}
D\left( M, d, \mu \right):=\sup_{x \in X, r > 0} \frac{\mu(B\left(x, 2r \right))}{\mu(B\left(x,r \right))} < \infty.
\end{equation}
The focus of this paper is the particular case where $\left( X, d, \mu
\right)=(K, d_{g}, \mu_{g})$ with the \textbf{volume doubling
  constant} denoted by $D_{g}:= D(K, d_{g}, \mu_{g})$.

Then, in the context of compact connected Lie groups, the two-sided spectral and heat kernel bounds in \eqref{e.1.0} and \eqref{e.1.1} would follow from the following conjecture.

\begin{con}\label{MC1} Let $K$ be a connected real compact Lie group. There is a constant $D(K)$
such that
\begin{equation}\label{e.1.4}
D_{g} \leqslant D(K) \text{ for all } g \in \mathfrak{L}(K),
\end{equation}
that is, $K$ is uniformly doubling with constant $D(K)$.
\end{con}

As an illustration of the significance of this conjecture,
the volume doubling constant also appears as the constant in the
Poincar\'e inequality (see Section
\ref{a.Poincare}):
 \begin{equation}\label{P1}
 \int_{B_{g}\left(x,r \right)}|f-f_{x, r}|^2 d\mu_{g} \leqslant   2 r^{2} D_{g} \int_{B_{g}\left(x,2r \right)} |\nabla_{g} f|_{g}^2 d\mu_{g} \text{ for all } f \in \mathcal C^\infty(B_{g}\left(x,2r \right)),
 \end{equation}
 where $f_{x, r}:=\int_{B_{g}\left(x,r \right)}f d\mu_{g}$ denotes the
 mean of $f$ over $B_{g}\left(x,r \right)$.  Hence, the validity of
 \eqref{e.1.4} implies that the constant in the Poincar\'e inequality
 \eqref{P1} is uniform over all metrics in $\mathfrak{L}(K)$.
 Together with known heat kernel estimates due to \cite{Grigoryan1991,
   Saloff-Coste1992b, Saloff-CosteBook2002} this shows that the
 validity of Conjecture \ref{MC1} implies that of the two-sided heat
 kernel bound \eqref{e.1.1}.  A simple test function argument shows
 that \eqref{e.1.4} also implies the spectral gap estimate in terms of
 the diameter as given in \eqref{e.1.0}.

In this article, we prove that Conjecture \ref{MC1} is valid for $K=\SU(2)$. Our main result is as follows.

\begin{theo} \label{thM}
There exists a constant $D$ such that, for any left-invariant Riemannian metric $g$ on $\SU(2)$, we have
$D_{g} \leqslant D$.
\end{theo}
Since the underlying manifold of $\SU(2)$ is the $3$-sphere $S^3$, this theorem provides uniform volume doubling for a large family of Riemannian metrics on $S^3$.  This holds despite the fact that the geometries $g \in \mathfrak{L}(\SU(2))$ are not uniformly bounded in other senses; for instance, even after rescaling to constant diameter, there is no universal lower bound for the Ricci curvatures of metrics $g \in \mathfrak{L}(\SU(2))$ as we discuss in Section \ref{s.Diameter}.

The proof of Theorem \ref{thM} is based on the following explicit description of the behavior of the volume growth function $V_{g}$. Each $g \in \mathfrak{L}(\SU(2))$ can be identified with an inner product on $\su(2)$; let $0< a_1 \leqslant a_2 \leqslant a_3<\infty$ denote the square roots of its three eigenvalues, with respect to the standard Euclidean structure on $\su(2)$ induced by the negative of the Killing form. We stress that these parameters depend on the metric $g$.

\begin{theo} \label{thV}
There are constants $b_1, b_2\in (0,\infty)$ such that for all $g \in \mathfrak{L}(\SU(2))$, the function $V_{g}$ satisfies
\[
b_1\leqslant \frac
{V_{g}\left( r \right)}
{\overline{V_{g}}\left( r \right)
}\leqslant b_2, \text{ where }
\]
\[
\overline{V_{g}}\left( r \right)=\left\{
\begin{array}{ll}
r^3 & \text{ if } 0 <r \leqslant a_1a_2/a_3,
\\
\left( a_3 / a_1a_2 \right)r^4 & \text{ if }  a_1a_2/a_3 < r\leqslant a_1,
\\
\left( a_1a_3 /a_2 \right) r^2 & \text{ if }  a_1 < r \leqslant a_2,
\\
a_1a_2a_3 & \text{ if } a_2< r <\infty.
  \end{array}\right.
  \]
\end{theo}
We note that $a_1$ can be characterized as the length of the shortest closed geodesic for $g$, while $a_2$ can be replaced in the theorem above by the diameter  $\operatorname{diam}_{g}(\SU(2))$ because the two are uniformly comparable (this is not entirely obvious, but will be proved in Section \ref{s.Diameter}), and
that $a_3$ is then uniformly comparable to the quantity $\mu_g(\SU(2))/a_1\operatorname{diam}_{g}$.

As far as we know, the only other case when Conjecture \ref{MC1} is
known to hold is for $K=\mathbb{T}^n$, the $n$-dimensional torus, for
any fixed $n$. This can be seen via lifting to the covering group,
$\mathbb{R}^n$, on which all Euclidean metrics are
isomorphic with the same doubling constant $2^n$.  But doubling passes
to quotients.  The key argument is given in \cite[Lemma
  1.1]{Guivarch1973a}; see also \cite[(5.5),
  p.20]{DiaconisSaloff-Coste1994}.  Alternatively, this can be seen
using curvature as explained in Section \ref{s.Curvature}, since every
left-invariant metric on a torus is flat and has zero Ricci curvature.

It is important to note that Theorem \ref{thM} implicitly includes two limit cases. In one case, the metric tends to infinity in one direction, and the manifold approaches  a sub-Riemannian manifold, which itself is doubling.  If the metric tends to zero in one direction, the 3-dimensional manifold $\SU(2)$ collapses to a 2-dimensional quotient, which is also doubling.  Then in some sense, the question becomes whether the doubling constant varies continuously with respect to these limits. One of the difficulties is that both cases must be considered simultaneously.

Our approach for $\SU(2)$ is rather explicit and makes use of its specific structure, with the important benefit of providing a detailed estimate of the volume function as stated in Theorem \ref{thV}.  We show that the volume function exhibits different behavior at different
scales: Euclidean behavior at very small scales, sub-Riemannian behavior at intermediate scales and ``quotient geometry'' behavior at relatively large scales, and this is done uniformly over all metrics in $\mathfrak{L}(K)$. This allows us to approximate the volume growth function of the metric $g$ by the simple explicit function $\overline{V_g}$ which essentially ``pieces together'' the growth functions of those three spaces.  We hope that the study of this special case will open the door to similar results for other compact groups.

\subsection{Curvature, or not}\label{s.Curvature}
In geometric analysis, ever since the pioneering work of S.-T.~Yau in the 1970s, Ricci curvature has been the tool of choice to prove spectral bounds and other analytic estimates such as various forms of Harnack inequalities and heat kernel estimates, especially if one is interested in statements that are uniform over large families of Riemannian manifolds. In particular, the celebrated Bishop-Gromov volume comparison theorem implies that for any complete Riemannian manifold $(M, g)$ of dimension at most $n$ with a non-negative Ricci curvature, the doubling constant $D(M, d_{g}, \mu_{g})$ is bounded by $2^n$, the doubling constant of Euclidean space $\mathbb{R}^{n}$.  If
the curvature condition is relaxed to a Ricci curvature lower bound, say, $\operatorname{Ric}_{g} \geqslant -\kappa g$, while keeping the restriction that the dimension is at most $n$, one still has a uniform bound on the doubling constant $D(M, d_{g},\mu_{g})$ as long as one
imposes a fixed upper bound on the diameter $\operatorname{diam}_{g}\left( M \right)$. In these contexts, the Poincar\'e inequality \eqref{P1} is not a direct consequence of the
doubling property, but it follows from the dimension and curvature assumptions (and an upper bound on the diameter in the case of $\operatorname{Ric}_{g}\geqslant -\kappa g$). In fact, fix  an $\epsilon >0$ and the dimension $n$.  For Riemannian manifolds of that fixed dimension, the curvature-diameter assumption
\[
\operatorname{Ric}_{g} \geqslant -\epsilon \operatorname{diam}_{g}^{-2}  g
\]
implies that $(M,g)$ is doubling and satisfies the Poincar\'{e} inequality (\ref{P1}) with constant depending only on $n$ and $\epsilon$. Note, however, that this curvature-diameter assumption is not invariant under multiplication  of the metric by a positive scalar. See the Bishop-Gromov comparison theorem and the result of P. Buser in \cite{Buser1982} and also \cite[Section 5.6.3]{Saloff-CosteBook2002}.

In this spirit, Conjecture \ref{MC1} is very much modeled on the non-negative Ricci curvature result
described above. Even so, except in the commutative case of the flat tori, it is  well known that no uniform Ricci lower bound  can hold over the entire family $\mathfrak{L}(K)$ of left-invariant metric on a group $K$.  In fact, the very nature of Conjecture \ref{MC1} implies that it not only covers left-invariant Riemannian geometries but also left-invariant sub-Riemannian geometries which can be described, in some rather obvious ways, as limits of left-invariant Riemannian geometries. This is made explicit for $\SU(2)$ in Section \ref{sub-riemannian}.

Recently there have been interesting attempts to extend curvature
techniques in the context of sub-Riemannian geometries
e.g. \cite{AgrachevBarilariRizzi2017, BaudoinGarofalo2017,
  BaloghTysonVecchi2016, CapognaCitti2016, Hladky2012}. However, even
in the case of left-invariant geometries on $\SU(2)$, it seems that
these curvature techniques (old and new) do not yield a proof of
Theorem \ref{thM}.

Other works have obtained geometric inequalities, including volume
doubling and the stronger measure contraction property $MCP(k,n)$
introduced by \cite{Ohta2007a}, that hold uniformly over a
one-parameter family of Riemannian geometries approximating a
sub-Riemannian geometry \cite{AgrachevLee2014,
  BaudoinGrongKuwadaThalmaier2017, Juillet2009, LeePaul2016a,
  LeePaulLiZelenko2016, Rifford2013}.  However, these works use very
different techniques, and all known results appear to rely on
assumptions of horizontal curvature bounds or additional symmetry,
such as Sasakian structure.  To the best of our knowledge, these
assumptions are not satisfied uniformly over all left-invariant
sub-Riemannian geometries on $\SU(2)$, and thus those results likewise do
not imply Theorem \ref{thM}.

\section{Preliminaries}\label{s.Prelim}

\subsection{The group $G=\SU(2)$ and left-invariant metrics on $G$}\label{s.2.1}

The compact Lie group $\SU(2)$ is the group of $2 \times 2$ complex matrices which are unitary and have determinant $1$.  The group identity of $\SU(2)$ is the identity matrix $I$, which we shall also denote by $e$ when emphasizing the group structure.  The corresponding Lie algebra $\su(2)$, identified with the tangent space $T_e \SU(2)$, is the space of $2 \times 2$ complex matrices which are skew-Hermitian and have trace $0$.  We note that a left-invariant metric $g$ on $\SU(2)$ is uniquely defined by its action on $\su(2)$, the tangent space at the identity.

Since $\SU(2)$ is compact, the Killing form $B(v,w) = \frac{1}{2} \tr(\operatorname{ad}_v \operatorname{ad}_w)$ is negative definite, and so $-B$ is an inner product on $\su(2)$ which is invariant.  It induces a bi-invariant Riemannian metric on $\SU(2)$, which we will
call the \textbf{canonical bi-invariant metric}; it is unique up to scaling because $\SU(2)$ is simple \cite[Lemma 7.6]{Milnor1976}.  In this canonical metric, $\SU(2)$ is isometric to a round
sphere.

As $\SU(2)$ is compact, by \cite[Lemma 7.2]{Milnor1976} a
left-invariant metric $g$ on $\SU(2)$ is bi-invariant if
and only if $\operatorname{ad}_{x}$ is skew-adjoint with respect to
$g$ for every $x \in \su(2)$. More detail (based mostly
on \cite{Milnor1976}) can be found in \cite[Chapter
  1.4]{ChowKnopfBook}.

\subsection{Standard Milnor bases}\label{s.MilnorBases}

A key property of $\SU(2)$ is that any left-invariant metric $g$ can
be diagonalized by a basis for $\su(2)$ for which the structure
constants have a very simple form.  Such bases were studied by Milnor
in \cite{Milnor1976}.

Throughout this section, $\left\{ i, j, k \right\}$ will be taken to
range over all cyclic permutations of the indices $\left\{ 1, 2, 3
\right\}$.

\begin{defin}\label{d.3.1}
  We shall say that an ordered basis $\left\{e_1, e_2, e_3 \right\}$ for $\su(2)$ is
  a \textbf{standard Milnor basis} if it satisfies the relations
  \begin{equation*}
    [e_1, e_2] = e_3, \qquad [e_2, e_3] = e_1, \qquad [e_3, e_1] = e_2,
  \end{equation*}
  or for short
  \begin{equation*}
    [e_i, e_j] = e_k.
  \end{equation*}
\end{defin}

\begin{exa}\label{pauli-ex}
  The Pauli matrices
  \begin{equation}
  \widehat{e}_1 = \frac{1}{2}
  \begin{pmatrix}
    0 & -i \\ -i & 0
  \end{pmatrix}, \quad
  \widehat{e}_2 = \frac{1}{2}
  \begin{pmatrix}
    0 & -1 \\ 1 & 0
  \end{pmatrix}, \quad
  \widehat{e}_3 = \frac{1}{2}
  \begin{pmatrix}
    -i & 0 \\ 0 & i
  \end{pmatrix},
\end{equation}
are a standard Milnor basis.
\end{exa}

\begin{exa}\label{milnor-transformations}
  If $\left\{e_1, e_2, e_3 \right\}$ is a standard Milnor basis, then so are
  \begin{enumerate}
  \item the cyclic permutations $\left\{ e_2, e_3, e_1\right\}$ and $\left\{ e_3, e_1,
    e_2 \right\}$;
  \item the ordered basis $\left\{-e_1, e_3, e_2 \right\}$.  As such, any
    permutation of a standard Milnor basis may itself be made into a
    standard Milnor basis by possibly negating one element;
  \item the basis
    \begin{equation}\label{milnor-rotate}
    \left\{ \cos(\theta) e_1 + \sin(\theta) e_2, -\sin(\theta) e_1 +
      \cos(\theta) e_2, e_3 \right\}, \qquad \theta \in \mathbb{R}.
    \end{equation}
  \end{enumerate}
\end{exa}

\begin{rem}\label{r.3.4} Definition \ref{d.3.1} is slightly different
  from a more common notion of Milnor frames, in which one begins with
  a metric $g$, and in addition to the commutation relations one
  assumes that $\left\{e_1, e_2, e_3 \right\}$ are orthogonal with
  respect to $g$.
\end{rem}

The next lemma is a consequence of the fact that all Lie algebra
automorphisms of $\su(2)$ are inner, and therefore the set of all
standard Milnor bases for $\su(2)$ coincides with the orbit of
$\operatorname{Ad}$ starting at any standard Milnor basis. Note that
this is not so for $\SU\left( n \right), n \geqslant
3$. As always for a matrix Lie group $G$ we use the fact that
$\operatorname{Ad}_{g} X=g Xg^{-1}$ for $g \in G$ and $X \in
\mathfrak{g}$, the Lie algebra of $G$, where on the right we have the
products of matrices.

\begin{lem}\label{automorphism}
  Suppose $\left\{ e_1, e_2, e_3 \right\}$  is a standard Milnor basis.  Then $\left\{ e_1^{\prime}, e_2^{\prime}, e_3^{\prime} \right\} \subseteq \su(2)$ is a standard Milnor basis if and only if there exists $y \in \SU(2)$ such that $\Ad_y e_i = y e_i y^{-1} = e_i^{\prime}$ for $i=1, 2, 3$.
\end{lem}

\begin{proof}
  For any $y \in \SU(2)$, the map $v \mapsto y v y^{-1}$ is a Lie algebra automorphism of $\su(2)$, so it is clear that $e_i^{\prime} = y e_i y^{-1}$ produces a standard Milnor basis.  Conversely, suppose $\left\{ e_1^{\prime}, e_2^{\prime}, e_3^{\prime} \right\}$ is a standard Milnor basis.  Since $\left\{ e_1, e_2, e_3 \right\}$ and $\left\{ e_1^{\prime}, e_2^{\prime}, e_3^{\prime} \right\}$ are both bases, there is a unique linear automorphism $T$ of the vector space $\su(2)$ satisfying $T e_i = e_i^{\prime}$, $i=1, 2, 3$.  Then if $(i, j, k)$ is any cyclic permutation
  of the indices $(1,2,3)$, we have
  \begin{equation*}
    [T e_i, T e_j] = [e_i', e_j'] = e_k' = T e_k = T[e_i, e_j].
  \end{equation*}
  It follows that $[Tu, Tv] = T[u,v]$ for any $u,v \in \{e_1, e_2, e_3 \}$, and by linearity the same holds for any $u,v \in \su(2)$.  So $T$ is a Lie algebra automorphism of $\su(2)$.  It is well-known that  every Lie algebra automorphism of $\su(2)$ is inner (i.e. the outer
  automorphism group is trivial) as pointed out in \cite[Proposition 5.1]{Wood1989a}. Thus $T = \Ad_y$  for some $y \in \SU(2)$.
\end{proof}

\begin{lem}\label{signs}
  Suppose $\left\{ e_1, e_2, e_3 \right\}$ is a basis for $\su(2)$ satisfying $[e_i,
    e_j] = \lambda_k e_k$ where $\lambda_i, \lambda_j, \lambda_k \in
  \{\pm 1\}$.  Then $\lambda_1 = \lambda_2 = \lambda_3$.  In
  particular, either $\left\{ e_1, e_2, e_3 \right\}$ or $\left\{ -e_1, e_2, e_3 \right\}$ is a
  standard Milnor basis.
\end{lem}

\begin{proof}
  Let $B(v,w) = \frac{1}{2} \tr(\operatorname{ad}_v
  \operatorname{ad}_w)$ be the Killing form of $\su(2)$, which is
  negative definite since $\SU(2)$ is compact.  Then a simple
  computation shows $B(e_i, e_i) = - \lambda_j \lambda_k$.  Since
  this must be negative for each $i$, it follows that $\lambda_1,
  \lambda_2, \lambda_3$ are all $+1$ or all $-1$.  In the former case,
  $\left\{ e_1, e_2, e_3 \right\}$ is already a standard Milnor basis,
  and in the latter case, it is easy to check that $\left\{ -e_1, e_2,
  e_3 \right\}$ is.
\end{proof}

\begin{lem}\label{masha-identities-lem}
  For any standard Milnor basis $\left\{ e_1, e_2, e_3 \right\}$, we have the following
  identities in the matrix algebra $M^{2 \times 2}(\mathbb{C})$
  \begin{equation}\label{masha-identities}
    e_i^2 = -\frac{1}{4} I, \qquad e_i e_j = \frac{1}{2} e_k, \qquad e_i e_j+e_j e_i=0,
  \end{equation}
  where $(i,j,k)$ is, as before, any cyclic permutation of the indices
  $(1,2,3)$ and $i \not= j$.
\end{lem}

\begin{proof}
Note that by Lemma \ref{automorphism}, it is enough to verify identities \eqref{masha-identities} for one standard Milnor basis since $\Ad_h I=I$ and $\Ad_h 0=0$ for all $h \in \SU(2)$. A simple calculation proves the first two identities for Pauli matrices, while the last identity can be shown by appealing to Definition \ref{d.3.1} and the second identity as follows
\begin{equation*}
  e_i e_j+e_j e_i=2e_i e_j-e_{k}=e_k-e_k=0.
\end{equation*}
\end{proof}

\subsection{Left-invariant Riemannian metrics on $\SU(2)$}

\begin{lem}\label{milnor-diagonalize}
  Let $g$ be any left-invariant metric on $\SU(2)$.  There exists a
  standard Milnor basis $\left\{ e_1, e_2, e_3 \right\}$ which is orthogonal in the metric $g$
  and satisfies $g(e_1, e_1) \leqslant g(e_2, e_2) \leqslant g(e_3, e_3)$.
\end{lem}

\begin{proof}
  Following \cite{Milnor1976}, we define a cross product $\times$ on the $3$-dimensional inner product space $(\su(2), g)$, unique up to a choice of orientation.  To see it another way, one can identify $(\su(2), g)$ with $(\mathbb{R}^3, \cdot)$, uniquely up to a choice  of orientation, and pull back the cross product from $\mathbb{R}^3$.
  As shown in \cite[Lemma 4.1]{Milnor1976}, there is a unique linear map $L$ on $\su(2)$ satisfying $L(u \times v) = [u,v]$, and it is self-adjoint with respect to $g$.  Let $\{ w_1, w_2, w_3 \}$ be a $g$-orthonormal basis of eigenvectors for $L$, with eigenvalues $\lambda_1, \lambda_2, \lambda_3$.  Reordering this basis if necessary, we can assume it is positively oriented, so that $w_i \times w_j = w_k$. Then
  \begin{equation*}
    [w_i, w_j] = L(w_i \times w_j) = L(w_k) = \lambda_k w_k.
  \end{equation*}
  Setting $e_i = |\lambda_j \lambda_k|^{-1/2} w_i$, we can verify that $[e_i, e_j] = \pm e_k$ for some choice of signs, and that $\{e_1, e_2, e_3\}$ is still $g$-orthogonal. Finally we can re-index this basis as needed so that $g(e_1, e_1) \leqslant g(e_2, e_2) \leqslant g(e_3, e_3)$.  By Lemma \ref{signs}, either $\left\{ e_1, e_2, e_3 \right\}$ or $\left\{ -e_1, e_2, e_3 \right\}$ is the  desired standard Milnor basis.
\end{proof}

\begin{notation}\label{n.MetricsParameters} For any left-invariant
  Riemannian metric $g$ on $\SU(2)$ let $a_1 \leqslant a_2 \leqslant
  a_3$ be the (ordered) square roots of the eigenvalues of the metric
  $g$ with respect to the canonical Euclidean form defined by the
  negative of the Killing form $B(v,w) = \frac{1}{2}
  \tr(\operatorname{ad}_v \operatorname{ad}_w)$.  We call
  $a_1,a_2,a_3$ the \textbf{parameters associated to the metric}
  $g$. For any $0 < a_1 \leqslant a_2 \leqslant a_3 < \infty$, let
  $g_{(a_1, a_2, a_3)}$ denote the unique left-invariant Riemannian
  metric on $\SU(2)$ for which
\begin{equation*}
    g_{(a_1, a_2, a_3)}(\widehat{e}_i, \widehat{e}_j) = a_i^2
    \delta_{ij}, \qquad i = 1,2,3,
\end{equation*}
where $\widehat{e}_i$ are the Pauli matrices defined in Example
\ref{pauli-ex}.  Since $B(e_i, e_j) = -\delta_{ij}$ for any standard
Milnor basis, the parameters associated to $g_{(a_1, a_2, a_3)}$ are
indeed $a_1, a_2, a_3$.  Note that $g_{(1,1,1)}$ is the canonical
bi-invariant metric.
\end{notation}

\begin{cor}\label{isometry}
   Let $g$ be any left-invariant metric on $\SU(2)$, and let $a_1,
   a_2, a_3$ be its parameters.  Then $(\SU(2), g)$ is isometrically isomorphic to  $(\SU(2), g_{(a_1, a_2,
    a_3)})$, where $g_{(a_1, a_2, a_3)}$ is as defined in Notation \ref{n.MetricsParameters}.
\end{cor}

\begin{proof}
  Choose a standard Milnor basis $\{ e_1, e_2, e_3 \}$
  which diagonalizes $g$ as in Lemma \ref{milnor-diagonalize}.  Since
  $\{e_1, e_2, e_3\}$ is orthonormal with respect to $-B$, we have
  $g(e_i, e_i) = a_i^2$.  The linear map $\varphi : \su(2) \to \su(2)$
  defined by $\varphi(e_i) = \widehat{e_i}$ is a Lie algebra
  automorphism, since both bases have the same structure constants.
  Since $\SU(2)$ is simply connected, $\varphi$ induces a Lie group
  automorphism of $\SU(2)$ whose differential at the identity is
  $\varphi$, which by construction is an isometry of the
  left-invariant metrics $g$ and $g_{(a_1, a_2, a_3)}$.
\end{proof}

\begin{rem} \label{ei-remark}
By Corollary \ref{isometry}, for each left-invariant Riemannian  metric with parameters $\left( a_1, a_2, a_3 \right)$,  there is a group isomorphism providing an isometry between $g_{(a_1, a_2, a_3)}$ and that metric. Hence it suffices to consider $g_{(a_1, a_2, a_3)}$. In what follows, we abuse notation and use $\left\{ e_{1}, e_{2}, e_{3} \right\}$  to denote both a general Milnor basis or the particular Milnor basis formed by the Pauli matrices.
\end{rem}

\subsection{Exponential identities}
Recall that we use $I$  for the identity matrix when we treat it as an element of the matrix space $M^{2\times 2}\left( \mathbb{C} \right)$. Whenever we want to emphasize the role of $I$ as the identity in the group $\SU(2)$ we use $e$.

\begin{lem}\label{A2}  For any $A \in \su(2)$, we have
\begin{equation*}
    A^2 =
    -\det(A) I.
    \end{equation*}
\end{lem}

\begin{proof}
  One can verify this by observing that a general matrix $A \in \su(2)$  is of the form $A = \left(
  \begin{smallmatrix}
    ai & b+ci \\ -b+ci & -ai
  \end{smallmatrix}
  \right)$, $a, b, c \in \mathbb{R}$ and computing directly.
\end{proof}

\begin{lem}\label{expA}
  For $A \in \su(2)$, we have
  \begin{equation}\label{expA-eq}
    \exp(A) = (\cos \rho) I + \frac{\sin \rho}{\rho} A,
  \end{equation}
  where $\rho^2 = \det A$.
\end{lem}

\begin{rem}\label{r.3.9}  First observe that this identity can be used also for $\rho=0$, since then $A=0$ and $\exp(A) = I$. This can be seen by using any standard Milnor basis and writing $A=ae_{1}+be_{2}+ce_{3}, a, b, c \in \mathbb{R}$. Then $\rho^{2}=\frac{1}{4}\left( a^{2}+ b^{2} + c^{2} \right)=\det A$. In particular, if $\rho=0$, then $a=b=c=0$.
\end{rem}

\begin{proof}
  Consider the expansion $\exp A = \sum_{n=0}^\infty \frac{A^n}{n!}$.
  Grouping even and odd terms we can write $\exp A = \sum_{k=0}^\infty
  \frac{A^{2k}}{(2k)!} + \sum_{k=0}^\infty \frac{A^{2k+1}}{(2k+1)!}$.  By Lemma \ref{A2}

\begin{align*}
 & A^{2k} = (-\rho^2 I)^k = (-1)^k \rho^{2k} I,
  \\
 & A^{2k+1} =
  (-1)^k \rho^{2k} A = \frac{(-1)^k \rho^{2k+1}}{\rho} A,
\end{align*}
so the first sum equals $(\cos \rho)I$ and the second equals $\frac{\sin \rho}{\rho} A$.
\end{proof}

\begin{lem}\label{log-formula}
  For any $x \in \SU(2) \setminus \{-I\}$, we have $x = \exp(A)$, where
  \begin{equation*}
    A = \frac{\rho}{\sin \rho}(x - (\cos \rho) I)
  \end{equation*}
  and
  \begin{equation*}
    \rho = \arccos\left(\frac{\tr x}{2}\right).
  \end{equation*}
\end{lem}

\begin{rem} Similarly to Remark \ref{r.3.9}, if $\rho=0$, so that $x=I$, we take $A=0$ which is consistent with this identity. For $\rho = \pi$ we have $x = -I$ and can take $A = 2 \pi e_1$, for instance.
\end{rem}

\begin{proof}
Let $A, \rho$ be as given.  Since $\cos \rho = \frac{\tr x}{2}$, it is apparent that $\tr A = 0$.  To see that $A$ is skew-Hermitian, note  that since $x$ is unitary with $\det x = 1$, Cayley-Hamilton gives
  \begin{equation*}
    x^* = x^{-1} = -x + (\tr x) I = -x + (2 \cos \rho) I.
  \end{equation*}
  As such,
  \begin{equation*}
    A + A^* = \frac{\rho}{\sin \rho} (x + x^* - (2 \cos \rho) I) = 0.
  \end{equation*}
  Hence $A \in \su(2)$.

  We now verify that $\det A = \rho^2$; then the result follows
  immediately from Lemma \ref{expA}.  Using Lemma \ref{A2} and the
  fact that $A^* = -A$ we have
  \begin{equation*}
    \det(A) I = A^* A = \frac{\rho^2}{\sin^2 \rho}((1 + \cos^2 \rho) I - \cos \rho(x + x^*))
  \end{equation*}
  since $x x^* = I$.  Taking traces and noting that $\tr x = \tr x^* =
  2 \cos \rho$, we have
  \begin{equation*}
    2 \det(A) = \frac{\rho^2}{\sin^2 \rho}(2 + 2 \cos^2 \rho - 4 \cos^2 \rho) = 2 \rho^2
  \end{equation*}
  as desired.
\end{proof}

\begin{lem}\label{masha-tech} Suppose $\left\{ e_{1}, e_{2}, e_{3} \right\}$ is a standard Milnor basis for $\su(2)$. Then
\[
e^{se_{1}}e^{te_{2}}e^{-se_{1}}=\exp\left( t\left(\cos s\right) e_{2}+
t\left(\sin s\right) e_{3} \right), s, t \in \mathbb{R}.
\]
\end{lem}

\begin{rem} The proof given below does not use $\su(2)$ specifically, only the commutation relations for the Milnor basis. For $\su(2)$ this result can also be shown directly by using Lemma \ref{expA} on both sides.
\end{rem}

\begin{proof} Let
\begin{align*}
& f\left( s \right):=e^{se_{1}} e_{2} e^{-se_{1}},
\\
& g\left( s \right):=\left(\cos s\right) e_{2}+ \left(\sin s\right) e_{3}.
\end{align*}
Then  we see that
\begin{align*}
& f^{\prime}\left( s \right)=e^{se_{1}}[e_{1}, e_{2}]e^{-se_{1}},
\\
& f^{\prime\prime}\left( s \right)=e^{se_{1}}[e_{1}, [e_{1}, e_{2}]]e^{-se_{1}}=-e^{se_{1}} e_{2} e^{-se_{1}}=-f\left( s \right);
\\
& f\left( 0 \right)=e_{2}, f^{\prime}\left( 0 \right)=[e_{1}, e_{2}]= e_{3};
\end{align*}
\begin{align*}
& g^{\prime}\left( s \right)=-\left(\sin s\right) e_{2}+ \left(\cos s\right) e_{3},
\\
& g^{\prime\prime}\left( s \right)=-\left(\cos s\right) e_{2} - \left(\sin s\right) e_{3}=-g\left( s \right);
\\
& g\left( 0 \right)=e_{2}, g^{\prime}\left( 0 \right)=e_{3},
\end{align*}
so by uniqueness of the initial value problem for ODEs these two functions coincide, that is,
\begin{equation}\label{e.3.1}
e^{s e_{1}} e_{2} e^{-se_{1}}=\left(\cos s\right) e_{2}+ \left(\sin s\right) e_{3}.
\end{equation}
Finally,
\begin{align*}
& e^{se_{1}}e^{te_{2}}e^{-se_{1}}=\sum_{n=0}^{\infty} \frac{t^{n}}{n!}e^{se_{1}}e_{2}^{n}e^{-se_{1}}
\\
& = \sum_{n=0}^{\infty} \frac{t^{n}}{n!}\left( e^{se_{1}}e_{2}e^{-se_{1}} \right)^{n}= \sum_{n=0}^{\infty} \frac{t^{n}}{n!}\left( \left(\cos s\right) e_{2}+ \left(\sin s\right) e_{3} \right)^{n}
\\
& =\exp\left(t \left( \left(\cos s\right) e_{2} + \left(\sin s\right) e_{3} \right) \right).
\end{align*}
\end{proof}

\begin{rem}
  By applying Lemma \ref{masha-tech} to the standard Milnor basis
  \begin{equation*}
    \{e_1,  \cos(\theta) e_2 + \sin(\theta) e_3, -\sin(\theta) e_2 + \cos(\theta) e_3\}
  \end{equation*}
  for $\theta \in \mathbb{R}$, as in \eqref{milnor-rotate}, we obtain
  the more general identity
  \begin{equation}\label{masha-tech-theta}
    \begin{split}
    \exp(s e_1) \exp(t (\cos(\theta) e_1 + \sin(\theta) e_2)) \exp(-s
    e_1) \\
    = \exp(t (\cos(\theta + s) e_2 + \sin(\theta + s) e_3)).
   \end{split}
  \end{equation}
\end{rem}

\begin{cor}\label{c.3.13} Let $\left\{ e_{1}, e_{2}, e_{3} \right\}$
  be a standard Milnor basis, and $A=\left( t\cos s \right) e_{2} +
  \left(t \sin s\right) e_{3}$. Letting $\rho^2 = \det A$ as in Lemma
  \ref{expA}, for such $A$ we have
\[
\rho=\frac{t}{2}
\]
as noted in Remark \ref{r.3.9}.  Then by Lemma \ref{expA}
\begin{align}
& e^{se_{1}}e^{te_{2}}e^{-se_{1}}= \exp A = \cos \rho I+ \frac{\sin \rho}{\rho}A \label{e.3.6}
\\
& =\left( \cos \frac{t}{2} \right)I+2\left( \sin \frac{t}{2} \cos s\right)e_{2}+2\left( \sin \frac{t}{2} \sin s\right)e_{3}. \notag
\end{align}
\end{cor}

\subsection{The volume function}\label{ss.3.5}
In what follows we take $0 < a_1 \leqslant a_2 \leqslant a_3 < \infty$. Recall that by Corollary \ref{isometry} it is enough to consider the left-invariant Riemannian metric $g_{(a_1, a_2, a_3)}$ on $\SU(2)$ defined in Notation \ref{n.MetricsParameters}.
\begin{notation} For the metric $g_{(a_1, a_2, a_3)}$ we denote by $d_{(a_1, a_2, a_3)}$  the corresponding Riemannian distance; by $B_{(a_1, a_2, a_3)}(x,r)$ we denote the open ball in the distance $d_{(a_1, a_2, a_3)}$ centered at $x$ of radius $r$; by $\mu_{(a_1, a_2, a_3)}$ we denote the Riemannian volume measure corresponding to $g_{(a_1, a_2, a_3)}$.
\end{notation}

\begin{notation}\label{n.ReferenceMeasure}
By $\mu_0$ we denote the bi-invariant Haar probability measure on $\SU(2)$.
\end{notation}
Then the  Riemannian volume measure $\mu_{(a_1, a_2, a_3)}$  is a constant multiple of $\mu_0$.  Specifically, we have $\mu_{(a_1, a_2, a_3)} = (16 \pi^2 a_1 a_2 a_3) \mu_0$. The constant can be found by observing that in the bi-invariant metric $g_{(1,1,1)}$, the group $\SU(2)$ is a round sphere whose circumference is $4 \pi$ as follows, for instance, from Lemma \ref{expA}.

\begin{notation}\label{Va-notation}
Let $V_{(a_1, a_2, a_3)}(r) = \mu_0(B_{(a_1, a_2, a_3)}(e, r))$ be the volume with respect to the measure $\mu_0$ of the ball in the metric $g_{(a_1, a_2, a_3)}$.
\end{notation}

Note that this is different from our previous notation $V_g$ used in Section \ref{s.MainResult}, since we are using the probability measure $\mu_0$ instead of the Riemannian volume measure $\mu_{(a_1, a_2, a_3)}$.  But this only makes a difference of a factor of $(16 \pi^2 a_1 a_2 a_3)^{-1}$, which for the purposes of studying volume doubling is irrelevant; and it is slightly more convenient for our purposes.

\begin{rem}\label{scaling-remark} For any $c > 0$, we have the scaling
\[
d_{(ca_1, ca_2, ca_3)}(x,y) = c d_{(a_1, a_2, a_3)}(x,y),
  \]
and so $B_{(ca_1, ca_2, ca_3)}(x,r) =   B_{(a_1, a_2, a_3)}(x, r/c)$.  As such, $g_{(a_1, a_2, a_3)}$ and  $g_{(ca_1, ca_2, ca_3)}$ have the same volume doubling constant. So for our purposes, we can suppose without loss of generality that $a_2 = 1$.  We show in Proposition \ref{diameter-bounds} that $a_2$ is comparable to the diameter of $g_{(a_1, a_2, a_3)}$, so the effect of this is rescaling of the metric to a roughly constant diameter.  The results in the remainder of the paper are written for general $a_2$, but in the proofs we generally work only with the case $a_2 = 1$.
\end{rem}

\begin{notation}\label{coords}
  Let $\Phi, \Psi : \mathbb{R}^3 \to \SU(2)$ be, respectively, the coordinates of the first and second kind (used in \cite{NagelSteinWainger1985}), defined by
  \begin{align*}
    \Phi(x_1, x_2, x_3) &= \exp\left(x_1 e_1 + x_2 e_2 + x_3 e_3\right),
    \\
    \Psi(y_1, y_2, y_3) &= \exp(y_1 e_1) \exp(y_2 e_2) \exp(y_3 e_3).
  \end{align*}
\end{notation}
We note that $\Phi, \Psi$ are both smooth maps, and that their differentials are isomorphisms at $(0,0,0)$.

\begin{notation}\label{J-notation}
Suppose $U \subset \mathbb{R}^3$ is open and $F : U \to \SU(2)$ is a diffeomorphism onto its image.  When we speak of the Jacobian $J : U \to (0,\infty)$ of $F$, we mean the normalization such that $\mu_0(F(K)) = \int_K J\,dm$ for measurable $K \subset U$. Here $m$ is the Lebesgue measure on $\mathbb{R}^{3}$.
\end{notation}

\begin{rem}\label{r.eta}
Let $\Phi, \Psi : \mathbb{R}^3 \to \SU(2)$ be coordinates of the first and second kind introduced in Notation \ref{coords}.  Since both $d\Phi(0,0,0)$ and $d\Psi(0,0,0)$ are nonsingular, then by the inverse function theorem, on some small box $(-\eta, \eta)^3$, both $\Phi$ and $\Psi$ are diffeomorphisms onto their images.  In particular, taking $\eta$ smaller if needed, their Jacobian determinants (with the normalization defined in Notation \ref{J-notation}) are bounded away from $0$ on $[-\eta, \eta]^3$. Therefore there is some universal constant $c$ such that for any measurable $K \subset (-\eta, \eta)^3$ we have
\begin{equation}\label{e.eta}
\mu_0(\Phi(K)) \geqslant c m(K), \hskip0.1in \mu_0(\Psi(K)) \geqslant c m(K),
\end{equation}
 where $m$ is the Lebesgue measure on $\mathbb{R}^3$ as before.
\end{rem}

\section{Euclidean regime}\label{s.EuclideanRegime}

At sufficiently small scales, the Riemannian manifold $(\SU(2),
g_{(a_1, a_2, a_3)})$ (with $0 < a_1 \leqslant a_2 \leqslant a_3 < \infty$) looks
like Euclidean space, so we expect the volume of a ball of radius $r$
to scale like $r^3$.  We need to determine, in terms of $a_1, a_2,
a_3$, how small the scale has to be to ensure this happens with a
uniform constant.

\begin{pro}\label{euclidean-pro}
  There are constants $c,C$ such that, uniformly in $a_1 \leqslant a_2 \leqslant a_3$, we have
  \begin{equation}
    c (a_1 a_2 a_3)^{-1} r^3 \leqslant V_{(a_1, a_2, a_3)}(r)
    \leqslant C (a_1 a_2 a_3)^{-1} r^3 \quad \text{ for }
    0 \leqslant r \leqslant a_1 a_2/a_3.
  \end{equation}
\end{pro}

An upper bound can be obtained from the form of the Bishop--Gromov
comparison theorem, and a direct computation of the Ricci curvature of
$g_{(a_1, a_2, a_3)}$.

\begin{lem}\label{bishop-gromov}
  Let $(M,g)$ be a $3$-dimensional complete Riemannian manifold with
  $\operatorname{Ric}_g \geqslant -\kappa g$.  Then for any $0 < s \leqslant r <
  \infty$ we have
  \begin{equation*}
    \frac{V_g(r)}{V_g(s)} \leqslant \left(\frac{r}{s}\right)^3 e^{\sqrt{2
        \kappa} r}.
  \end{equation*}
\end{lem}

\begin{proof}
  By the Bishop--Gromov comparison theorem (see \cite[Corollary 5.6]{EschenburgBook} or \cite[Lemma 36]{PetersenBook2ndEdition}), we have
  \begin{equation}
    \frac{V_g(r)}{V_g(s)} \leqslant \frac{V_{\kappa}(r)}{V_{\kappa}(s)},
  \end{equation}
where $V_{\kappa}(r)$ is the volume of a ball of radius $r$ in the 3-dimensional hyperbolic space of constant sectional curvature  $-\kappa/2$ (which has constant Ricci curvature $-\kappa$). The
  volume $V_{\kappa}(r)$ is given by \cite{Wielenberg1981}
  \begin{equation*}
   V_{\kappa}(r)= \pi (\kappa/2)^{-3/2} \left(\sinh(\sqrt{2 \kappa} r) - \sqrt{2 \kappa} r\right)
  \end{equation*}
and the desired result follows by observing that
  \begin{equation*}
    \frac{x^3}{6} \leqslant \sinh(x) - x \leqslant \frac{x^3 e^x}{6}, \qquad x
    \geqslant 0
  \end{equation*}
 which can be seen, for instance, by inspecting the Taylor series.
\end{proof}

\begin{proof}[Proof of Proposition \ref{euclidean-pro}]
It is enough to bound the Ricci tensor of the metric $g_{(a_1, a_2, a_3)}$.  In the basis $\{e_1, e_2, e_3\}$, $\Ric$ is diagonal, and we find
\[
\Ric(e_i,e_i) = \frac{\left( a_i^4 -(a_j^2-a_k^2)^2 \right)}{2(a_j a_k)^{2}},
\]
where $(i,j,k)$ is any permutation of $(1,2,3)$ (note that the
expression is symmetric in $a_j$ and $a_k$, so it is not necessary to
restrict to positive permutations). Now we need to find the smallest of the ratios
\[
\frac{\Ric (e_i,e_i)}{g(e_i,e_i)}=
\frac{\left( a_i^4 -(a_j^2-a_k^2)^2\right)}{2(a_ia_j a_k)^{2}}.
\]
Recall that $a_1 \leqslant a_2 \leqslant a_3$, and therefore $|a_j^2 - a_k^2|
\leqslant a_3^2$, and so we have  $a_i^4 -(a_j^2-a_k^2)^2 \geqslant -a_3^4$.  This yields the bound
\begin{equation}\label{ricci-bound}
  \frac{\Ric (e_i,e_i)}{g(e_i,e_i)} \geqslant -\frac{a_3^4}{2(a_i a_j a_k)^{2}} = -\frac{1}{2}
  \left(\frac{a_3}{a_1 a_2}\right)^2
\end{equation}
which is sharp when $i=3$ and $a_1 = a_2$.  Let us denote by
$\kappa:= \frac{1}{2} \left(\frac{a_3}{a_1 a_2}\right)^2$ the quantity on the
right side of \eqref{ricci-bound}.  If $r \leqslant \frac{a_1
  a_2}{a_3}$, then we have $\sqrt{2 \kappa} r \leqslant 1$, and Lemma
\ref{bishop-gromov} gives
\begin{equation} \label{Vg-euclidean-bound}
  V_{g}(r) \leqslant e r^3 s^{-3} V_{g}(s), \qquad 0 < s \leqslant r \leqslant \frac{a_1 a_2}{a_3}.
\end{equation}
Letting $s \to 0$, we have $V_g(s) \sim \frac{4}{3} \pi s^3$ (since a Riemannian manifold is locally Euclidean), so that \eqref{Vg-euclidean-bound} reads $V_g \leqslant C r^3$ where $C = \frac{4}{3} \pi e$. Rewriting this in terms of $V_{(a_1, a_2, a_3)}$ using Notation \ref{Va-notation}, we have the upper bound
\begin{equation*}
    V_{(a_1, a_2, a_3)}(r) \leqslant C (a_1 a_2 a_3)^{-1} r^3, \qquad 0 \leqslant r \leqslant \frac{a_1 a_2}{a_3}
\end{equation*}
absorbing $1/16\pi^2$ into the constant $C$.

Now we turn to the lower bound.
Let $\Psi : \mathbb{R}^3 \to \SU(2)$ be coordinates of the second kind introduced in Notation \ref{coords}. By Remark \ref{r.eta}, there exist $\eta > 0$ and a constant $c$ such that for any measurable $K \subset (-\eta, \eta)^3$ we have $\mu_0(\Psi(K)) \geqslant c m(K)$, where $m$ is the Lebesgue measure on $\mathbb{R}^3$. Suppose that $t \leqslant \eta a_1$ and consider the box
\begin{equation*}
  K_t = \pminterval{\frac{t}{a_1}} \times
  \pminterval{\frac{t}{a_2}} \times
  \pminterval{\frac{t}{a_3}} \subset
  [-\eta, \eta]^3.
\end{equation*}
On the one hand, we have $\mu_0(\Psi(K_t)) \geqslant c m(K_t) = 8 c
(a_1 a_2 a_3)^{-1} t^3$.  On the other hand, for any $(x,y,z) \in
K_t$, we have
\begin{equation*}
  d_{(a_1, a_2, a_3)}(e, \Psi(x,y,z)) \leqslant a_1 |x| + a_2 |y| + a_3 |z|
  \leqslant 3 t.
\end{equation*}
That is, $\Psi(K_t) \subset B_{(a_1, a_2, a_3)}(3 t)$, so we have
\begin{equation*}
  V_{(a_1, a_2, a_3)}(3t) = \mu_0(B_{(a_1, a_2, a_3)}(3t)) \geqslant   \mu_0(\Psi(K_t)) \geqslant 8c(a_1 a_2 a_3)^{-1} t^3
\end{equation*}
or, letting $r = 3t$,
\begin{equation}
  V_{(a_1, a_2, a_3)}(r) \geqslant c^{\prime} (a_1 a_2 a_3)^{-1} r^3, \quad 0
  \leqslant r \leqslant\frac{\eta}{3}a_1,
\end{equation}
where $c^{\prime} = \frac{8}{27} c$.  To complete the proof for all $0
\leqslant r \leqslant a_1$, note that for any $\frac{\eta}{3} a_1
\leqslant r \leqslant a_1$ we have by the monotonicity of $V$ that
\begin{equation*}
  V_{(a_1, a_2, a_3)}(r) \geqslant V_{(a_1, a_2, a_3)}\left(\frac{\eta}{3} a_1\right)
  \geqslant c' \frac{\eta^3}{27} \frac{a_1^2}{a_2 a_3} \geqslant c^{\prime}
  \frac{\eta^3}{27} (a_1 a_2 a_3)^{-1} r^3
\end{equation*}
so taking $c^{\prime \prime} = \min\{1, \frac{\eta^3}{27}\} c^{\prime}$ we have the desired
\begin{equation}
  V_{(a_1, a_2, a_3)}(r) \geqslant c^{\prime \prime} (a_1 a_2 a_3)^{-1} r^3, \quad 0
  \leqslant r \leqslant  a_1
\end{equation}
and in particular for $0 \leqslant r \leqslant \frac{a_1 a_2}{a_3}$, since $a_2
\leqslant a_3$.
\end{proof}

\section{Heisenberg regime}\label{s.HeisenbergRegime}

For $r \geqslant a_1 a_2/a_3$, the Euclidean behavior breaks down.
The growth of a ball in the $e_3$ direction is now affected by the
relation $[e_1, e_2] = e_3$; paths can make more efficient progress in
the $e_3$ direction by making a loop in the $e_1$ and $e_2$
directions.  This is well approximated by the sub-Riemannian geometry
of the $3$-dimensional Heisenberg group, in which one cannot move
tangent to the vertical direction $e_3$ at all.  The sub-Riemannian
Heisenberg group has Hausdorff dimension $4$, which accounts for the
$r^4$ volume scaling that appears in this regime.


\begin{lem}\label{H-lemma}
  Define $H : \mathbb{R}^2 \to \SU(2)$ by
  \begin{equation}\label{H-def}
    H(u,v) = \exp(-u e_1) \exp(-v e_2) \exp(u e_1) \exp(v e_2).
  \end{equation}
Then in some neighborhood $U$ of $(0,0)$ in $\mathbb{R}^2$ we can write
\begin{equation}
  H(u,v) = \exp(uv h(u,v))
\end{equation}
where $h : \mathbb{R}^2 \to \su(2)$ is $C^\infty$ and
$h(0,0) = e_3$.
\end{lem}

\begin{proof} We give two different arguments. Applying the Campbell--Baker--Dynkin--Hausdorff formula gives a power series for $\log H(u,v)$, convergent in a neighborhood $U$ of $(0,0)$.  The first-order terms in $u,v$ cancel out, while the second-order  term is $uv e_3$.  Each higher-order term consists of a combinatorial coefficient multiplied by some iterated Lie bracket of the vectors $u e_1$ and $v e_2$.  In any term that does not vanish, the innermost bracket must be of the form $\pm [u e_1, v e_2] = \pm uv e_3$; so such an iterated bracket must equal $\pm u^a v^b e_i$ with $a,b \geqslant 1$.  Thus we can factor out $uv$ from every term of this power series, writing $\log H(u,v) = uv h(u,v)$ where $h$ is given by a convergent power series and thus is real analytic in $U$.

Alternatively, one can give a more direct proof by using \eqref{e.3.6}, \eqref{expA-eq}, \eqref{masha-identities} and the double angle formula to write
  \begin{equation}\label{H-basis-formula}
    \begin{split}
      & H\left( u, v \right)= \frac{1}{2}\left( \left( 1 +\cos u\right)  + \cos v \left( 1 -\cos u\right) \right)I
      -\left( \left( 1 - \cos v \right) \sin u\right)e_{1}
      \\
      &
      + \sin v\left(1-\cos u \right) e_{2}
      +\left( \sin v \sin u\right)e_{3}.
  \end{split}
  \end{equation}
Then by using Lemma \ref{log-formula}, one obtains a formula for $\log(H(u,v))$.  It can then be seen by inspection that $h(u,v) = \log(H(u,v))/(uv)$ has a removable singularity at $(0,0)$, where the limit equals $e_3$.
\end{proof}


\begin{lem}\label{F-lemma}
  Define $F : \mathbb{R}^4 \to \SU(2)$ by
  \begin{equation*}
    F(s_1, s_2, s_3,\delta) = \exp(s_1 e_1) \exp(s_2 e_2) H\left(\sgn(s_3)
    \sqrt{|s_3|}, \delta \sqrt{|s_3|}\right).
  \end{equation*}
Then there is a neighborhood $V$ of $(0,0,0) \in \mathbb{R}^3$ such that on $V \times [0,1]$, the partial derivatives of $F$ with respect to $s_1, s_2, s_3$ exist and are jointly continuous, and we  have
  \begin{equation*}
    \partial_{s_1} F(0,0,0,\delta) = e_1, \qquad \partial_{s_2}
    F(0,0,0,\delta) = e_2, \qquad \partial_{s_3} F(0,0,0,\delta) =
    \delta e_3.
  \end{equation*}
Moreover, there is a jointly continuous $f : V \times [0,1] \to T \SU(2)$ such that
  \begin{equation*}
    \partial_{s_3} F(s_1, s_2, s_3, \delta) = \delta
    f(s_1, s_2, s_3, \delta)
  \end{equation*}
and $f(0,0,0,\delta) = e_3$.
\end{lem}

\begin{proof}
Observe that $s \mapsto \sgn(s) \sqrt{|s|}$ is continuous on $\mathbb{R}$, and $H$ is smooth, with $H(0,0) = I$.  Thus the desired statements about $\partial_{s_1} F$, $\partial_{s_2} F$ are clear.

For convenience, let $G(s,\delta) = H\left(\sgn(s) \sqrt{|s|}, \delta \sqrt{|s|}\right)$, so that $F(s_1, s_2, s_3, \delta) = \exp(s_1 e_1) \exp(s_2 e_2) G(s_3, \delta)$.  Let us write $H(u,v) = \exp(uv \,h(u,v))$ as in the previous lemma.  Then for $s$ in some interval $(-\epsilon, \epsilon)$ we can write
  \begin{equation*}
    G(s, \delta) = \exp(\delta k(s,\delta))
  \end{equation*}
  where
  \begin{equation*}
    k(s, \delta) = s \cdot h\left(\sgn(s)
    \sqrt{|s|}, \delta \sqrt{|s|}\right).
  \end{equation*}
  For $s \ne 0$, we compute
  \begin{align*}
    \partial_{s} k(s, \delta) &= h\left(\sgn(s)
    \sqrt{|s|}, \delta \sqrt{|s|}\right) \\
    &\quad + \frac{1}{2} \sgn(s) \sqrt{|s|} \,\partial_u h\left(\sgn(s)
    \sqrt{|s|}, \delta \sqrt{|s|}\right) \\
    &\quad + \frac{\delta}{2} \sqrt{|s|}  \,\partial_v h\left(\sgn(s)
    \sqrt{|s|}, \delta \sqrt{|s|}\right).
  \end{align*}
As $s \to 0$, the right side approaches $h(0,0) = e_3$, uniformly in $\delta \in [0,1]$.  Since $k$ is continuous, it follows (by L'H\^opital's rule) that $\partial_s k(0, \delta)$ exists and equals $e_3$; moreover, $\partial_s k$ is jointly continuous on $(-\epsilon, \epsilon) \times [0,1]$.

Now from the chain rule, since $\exp$ is smooth, we conclude that $\partial_s G(s,\delta)$ exists on $(-\epsilon, \epsilon) \times [0,1]$ and is given by
  \begin{equation*}
    \partial_s G(s,\delta) = d \exp_{\delta k(s,\delta)} \left[
    \partial_s [\delta k(s, \delta)] \right] = \delta \cdot d \exp_{\delta k(s,\delta)} \left[
    \partial_s k(s, \delta) \right]
  \end{equation*}
where $d \exp_{\delta k(s,\delta)} \left[ \partial_s k(s, \delta) \right]$ is a jointly continuous function of $s$ and $\delta$.  It is also clear from this that $\partial_s G(0,\delta) = \delta e_3$.  The desired statements about $\partial_{s_3} F$ follow.
\end{proof}

\begin{lem}\label{J-lemma}
For $\delta \in [0,1]$, consider $F^\delta = F(\cdot, \cdot, \cdot, \delta)$ as a map from $\mathbb{R}^3$ to $\SU(2)$.  Let $J^\delta$ be its Jacobian determinant as in Notation  \ref{J-notation}. Then there is a neighborhood $W$ of $(0,0,0) \in \mathbb{R}^3$ and a constant $c > 0$, independent of $\delta$, such that $J^\delta \geqslant c \delta$ on $W$.  In particular, for any measurable $K \subset W$, we have $\mu_0(F^\delta(K)) \geqslant c \delta m(K)$.
\end{lem}

\begin{proof}
Let $\omega$ be the Riemannian volume form on $\SU(2)$ associated to the bi-invariant metric $g_{(1, 1, 1)}$. Then we have
  \begin{align*}
    J^\delta &= \frac{1}{16 \pi^2} \omega(\partial_{s_1} F^\delta,
    \partial_{s_2} F^\delta, \partial_{s_3} F^\delta).
  \end{align*}
  If we set
  \begin{equation*}
    j(s_1, s_2, s_3, \delta) = \omega(\partial_{s_1} F^\delta(s_1,
    s_2, s_3),
    \partial_{s_2} F^\delta(s_1, s_2, s_3), f(s_1, s_2, s_3, \delta))
  \end{equation*}
where $f$ is as in Lemma \ref{F-lemma}, then $J^\delta = \frac{\delta}{16 \pi^2} j$.  Moreover, $j$ is jointly continuous on $V \times [0,1]$, and we have $j(0,0,0,\delta) = \omega(e_1, e_2, e_3) = 1$ for all $\delta$.  As such, by continuity and the compactness of $[0,1]$, there is a neighborhood $W \subset V$ of $(0,0,0) \in \mathbb{R}^3$ such that $j \geqslant \frac{1}{2}$ on $W \times [0,1]$, which implies $J \geqslant \frac{1}{32 \pi^2} \delta$.
\end{proof}

\begin{pro}[Heisenberg type lower bound]\label{heisenberg-lower} There is a constant $c$ such that, uniformly in $a_1\leqslant a_2\leqslant a_3$,
\[
 V_{(a_1, a_2, a_3)}(r) \geqslant c (a_1a_2)^{-2} r^4  \quad \text{ for }  0\leqslant r \leqslant  a_1.
\]
\end{pro}
Note that for $r\simeq a_1 a_2/a_3$ this lower bound matches the result
provided by Proposition \ref{euclidean-pro}.

\begin{proof}
Since the right side is consistent with the scaling described in Remark \ref{scaling-remark}, we suppose without loss of generality that $a_2 =
  1$.

  Let $F^\delta$ be as in Lemma \ref{F-lemma} and $W$ as in Lemma
  \ref{J-lemma}.  Choose $\eta > 0$ so small that $[-\eta, \eta]^3
  \subset W$.  We note that
  \begin{equation*}
    d_{(a_1, 1, a_3)}(e, F^\delta(s_1, s_2, s_3)) \leqslant s_1 a_1 + s_2 +
    2 a_1 \sqrt{|s_3|} + 2 \delta \sqrt{|s_3|}.
  \end{equation*}
  Now let us take $\delta = a_1 \in [0,1]$, so that this becomes
  \begin{equation*}
    d_{(a_1, 1, a_3)}(e, F^\delta(s_1, s_2, s_3)) \leqslant s_1 a_1 + s_2 +
    4 a_1 \sqrt{|s_3|}.
  \end{equation*}
  Suppose $r \leqslant a_1 \eta \leqslant \eta$ and let
  \begin{equation*}
    K_r = \pminterval{\frac{r}{a_1}}
    \times \pminterval{r} \times \pminterval{\frac{r^2}{a_1^2}}
  \end{equation*}
  so that $K_r \subset [-\eta, \eta]^3 \subset W$. We then have $m(K_r) = 8 a_1^{-3} r^4$ and $F^{a_1}(K_r) \subset B_{(a_1, 1, a_3)}(6r)$.  By Lemma \ref{J-lemma}, we have
  \begin{align*}
    V_{(a_1, 1, a_3)}(6r) = \mu_0(B_{(a_1, 1, a_3)}(6r)) \geqslant \mu_0(F^{a_1}(K_r)) \geqslant c a_1 m(K_r) = 8c
  a_1^{-2} r^4.
  \end{align*}
  or
  \begin{equation*}
    V_{(a_1, 1, a_3)}(r) \geqslant c' a_1^{-2} r^4, \quad 0 \leqslant r \leqslant 6 \eta a_1
  \end{equation*}
  where $c' = 8 c / 6^4$.  If it happens that $6 \eta \geqslant 1$ then we
  are finished; if not, we can drop the $6 \eta$ in the
  upper limit on $r$ as in the proof of Proposition
  \ref{euclidean-pro}, replacing $c'$ by  $(6 \eta)^4 c'$.
\end{proof}

\begin{pro}[Heisenberg type upper bound]\label{heisenberg-upper} There exists $\eta\in (0,1)$
 and a constant $C < \infty$ such that, uniformly in $a_1\leqslant a_2\leqslant a_3$,
\[
 V_{(a_1, a_2, a_3)}(r) \leqslant C \left( (a_1a_2a_3)^{-1}r^3 +  (a_1a_2)^{-2} r^4\right)
 \quad\text{ for } 0\leqslant r\leqslant \eta a_1.
 \]
 In particular, we have
\[
 V_{(a_1, a_2, a_3)}(r) \leqslant 2C (a_1a_2)^{-2} r^4 \quad\text{ for } a_1 a_2/a_3 \leqslant r\leqslant \eta a_1.
 \]
\end{pro}

\begin{proof}  Again, we assume $a_2 = 1$. Suppose $r \leqslant \eta a_1$, where $\eta$ is to be chosen later, and let $g \in B_{(a_1, 1, a_3)}(r)$.  This means that there is a smooth path $\gamma : [0,1] \to \SU(2)$ with $\gamma(0) = e$, $\gamma(1) =
  g$, and length $\ell_{(a_1, 1, a_3)}[\gamma] < r$.  Reparametrizing
  $\gamma$ by constant speed (with respect to $g_{(a_1, 1, a_3)}$), we
  can write $\dot{\gamma}(t) = \sum_{i=1}^3 \lambda_i(t)
  \widetilde{e_i}(\gamma(t))$, where $\widetilde{e_i}$ is the
  left-invariant vector field which equals $e_i$ at the identity, and
  $\sum_{i=1}^3 |a_i \lambda_i(t)|^2 \leqslant r^2$ for all $t \in [0,1]$.
  In particular, $|\lambda_i(t)| \leqslant r/a_i$.

We now invoke a result of R.~Strichartz \cite{Strichartz1987a} which extends the Baker--Campbell--Hausdorff--Dynkin formula by giving an exact expression for the exponential coordinates of $g$ in terms of $\lambda_i$.  The Strichartz (or Chen--Strichartz) formula says that $g = \exp z$, where
\begin{equation}\label{bcdhs}
    z = \sum_{n=1}^\infty \sum_{I\in \{1,2,3\}^n}
\left(\sum_{\sigma\in S_n}
\left(\frac{(-1)^{e(\sigma)}}{n^2 \binom{n-1}{e(\sigma)}}\right)
\int_{\Delta^n} \prod_{m=1}^n\lambda_{i_m}(s_{\sigma(m)})ds
\right) e_I \in \su(2).
  \end{equation}
  Here $I = (i_1, \dots, i_n)$, and $e_I$ is the $n$-fold iterated
  bracket
\begin{equation*}
e_I = [[\dots[e_{i_1}, e_{i_2}], \dots ], e_{i_n}].
\end{equation*}
Note that since $\{ e_i \}$ is a standard Milnor basis, each $e_I$ equals either 0 or some $\pm e_k$.  $S_n$ is the set of permutations of $\{1, \dots, n\}$, and  following Strichartz's notation, $e(\sigma) = |\{m < n : \sigma(m+1) < \sigma(m)\}|$ denotes the number of ``errors'' (also called ``descents'') of the permutation $\sigma$; for our purposes, we need only note that $e(\sigma)$ is an integer between $0$ and $n-1$.  Finally, $\Delta^n \subset [0,1]^n$ is the standard $n$-simplex $\{0 \leqslant s_1 \leqslant \dots \leqslant s_n \leqslant 1\}$, whose volume is $1/n!$.

Let us write $z = \sum_{i=1}^3 z_i e_i$.  We shall bound each of the $|z_i|$, which will show that $g$ is contained in the image under the coordinates $\Phi$ (see Notation \ref{coords}) of some box in $\mathbb{R}^3$ of bounded size.  This fact, combined with Remark \ref{r.eta} on the Jacobian determinant of $\Phi$, will give us an upper volume estimate for $B_{(a_1, 1, a_3)}(r)$.

We begin with $z_1$; the analysis of $z_2, z_3$ will be similar. Let $\zeta_{i,n}$ be the coefficient of $e_i$ in the $n$ term of the sum in \eqref{bcdhs}, so that $z_1 = \sum_{n=1}^\infty \zeta_{1,n}$. We must consider which values of $I$ give $e_I = \pm e_1$.  For $n=1$ we have only $I=(1)$, and for $n=2$ we have $I=(2,3)$ and $I=(3,2)$.  So we have
  \begin{align*}
    \zeta_{1,1} &= \int_0^1 \lambda_1(s)\,ds \\
    \zeta_{1,2} &= \frac{1}{4} \int_{0 \leqslant s_1 \leqslant
      s_2 \leqslant 1} (\lambda_2(s_1) \lambda_3(s_2) - \lambda_3(s_1)
    \lambda_2(s_2))\,ds_1 \,ds_2.
  \end{align*}
  This trivially gives
  \begin{equation}\label{zeta11-12}
    |\zeta_{1,1}| \leqslant \frac{r}{a_1}, \qquad |\zeta_{1,2}|
  \leqslant \frac{1}{4} \frac{r^2}{a_3}.
  \end{equation}

For $n \geqslant 3$, in order to have $e_I = \pm e_1$ we note that $i_1, i_2$ cannot both equal $1$ (else $e_I = 0$), and $i_n$ cannot equal $1$ either (since $[e_k, e_1] \ne \pm e_1$ for any $k=1,2,3$). So at least two of the $i_m$ are different from $1$, meaning that the corresponding $\lambda_{i_m}$ are bounded by $r$.  Since $|\lambda_i| \leqslant r/a_i$ and $a_1 \leqslant a_2 = 1 \leqslant a_3$, the remaining $\lambda_{i_m}$ are bounded by $r/a_1$, and we conclude that $\left|\prod_{m=1}^n \lambda_{i_m}(s_{\sigma(m)})\right| \leqslant r^n/a_1^{n-2}$.

Now to estimate the value of the parenthesized sum over $\sigma \in S_n$ in \eqref{bcdhs}, we note that $\Delta^n$ has a volume of $1/n!$, that $|S_n| = n!$, and that the combinatorial coefficient is at most 1. So this sum is bounded by $r^n/a_1^{n-2}$ as well. Finally, the total number of $I \in \{1,2,3\}^n$ is $3^n$, even though most of these do not yield $e_I = \pm e_1$.  So we conclude
  \begin{equation}
    |\zeta_{1,n}| \leqslant 3^n \frac{ r^n}{a_1^{n-2}} = 9 r^2 \left(3
    \frac{r}{a_1}\right)^{n-2}, \qquad n \geqslant 3.
  \end{equation}
By taking, say, $\eta < \frac{1}{6}$, so that $3 \frac{r}{a_1} <
\frac{1}{2}$, we can conclude
  \begin{equation}\label{zeta1n}
  \sum_{n=3}^\infty |\zeta_{1,n}|
  < 9 r^2 \sum_{n=3}^\infty \left(\frac{1}{2}\right)^{n-2} = 9r^2.
  \end{equation}
  Combining \eqref{zeta11-12} and \eqref{zeta1n}, we see that
  $\frac{r}{a_1}$ dominates, and we have
  \begin{equation}
    |z_1| \leqslant \sum_{n=1}^\infty |\zeta_{1,n}| \leqslant c \frac{r}{a_1}
  \end{equation}
for some universal constant $c$ ($c=11$ would do).

By similar arguments, we can obtain
  \begin{equation}
    |\zeta_{2,1}| \leqslant r,  \quad |\zeta_{2,2}| \leqslant
    \frac{1}{4} \frac{r^2}{a_1 a_3}.
  \end{equation}
Since $\frac{r}{a_1} \leqslant \eta \leqslant 1$ and $a_3 \geqslant 1$, both terms are dominated by $r$.  To estimate $\zeta_{2,n}$ for $n \geqslant 3$, we use the cruder fact that in order to get $e_I \ne 0$, we must have either $i_1$ or $i_2$ different from 1.  This leads to the estimate
  \begin{equation*}
    |\zeta_{2,n}| \leqslant 3^n \frac{r^n}{a_1^{n-1}} = 9 \frac{r^2}{a_1} \left(3
    \frac{r}{a_1}\right)^{n-2}
  \end{equation*}
and thus, stil with $\eta < \frac{1}{6}$,
  \begin{equation*}
    \sum_{n=3}^\infty |\zeta_{2,n}| < 9 \frac{r^2}{a_1}
  \end{equation*}
which again is dominated by $r$.  So we have shown
  \begin{equation}
    |z_2| \leqslant c r
  \end{equation}
increasing the value of the universal constant $c$ as needed.

For $z_3$, we obtain
  \begin{equation*}
    |\zeta_{3,1}| \leqslant \frac{r}{a_3}, \quad |\zeta_{3,2}| \leqslant
    \frac{1}{4} \frac{r^2}{a_1}
  \end{equation*}
and as before
  \begin{equation*}
    \sum_{n=3}^\infty |\zeta_{3,n}| < 9 \frac{r^2}{a_1}.
  \end{equation*}
We conclude
  \begin{equation}
    |z_3| \leqslant c \left(\frac{r}{a_3} + \frac{r^2}{a_1}\right),
  \end{equation}
where the first term dominates when $r \ll a_1/a_3$.

  As such, if we let
  \begin{equation*}
    K_r = \pminterval{c \frac{r}{a_1}} \times
    \pminterval{cr} \times \pminterval{c \left(\frac{r}{a_3} +
      \frac{r^2}{a_1}\right)}
  \end{equation*}
  so that in particular we have $K_r \subset [-c\eta, c\eta]^3$, we have
  that $B_{(a_1, 1, a_3)} \subset \Phi(K_r)$.  Letting $M$ be the
  maximum of the Jacobian determinant of $\Phi$ over $[-c\eta,
    c\eta]^3$, we have
  \begin{equation*}
    V_{(a_1, 1, a_3)}(r) \leqslant M m(K_r) = 8 M c^3 \left(\frac{r^3}{a_1
      a_3} + \frac{r^4}{a_1^2}\right)
  \end{equation*}
  which is the desired bound.
\end{proof}

\section{After Heisenberg}\label{s.CollapseRegime}

When $r$ exceeds $a_1$, the global geometry of $\SU(2)$ becomes
important.  Our ``budget'' $r$ is now large enough to let us travel
all the way around the sphere $\SU(2) \cong S^3$ in the ``cheap''
$e_1$ direction, and nothing is gained by traveling around the sphere
more than once.  So travel in the $e_1$ direction has negligible cost,
and the volume growth is comparable to what happens if we actually set
$a_1 = 0$.  The group $\SU(2)$ would collapse to a coset space mod the
subgroup $S = \{\exp(s e_1) : s \in \mathbb{R}\}$ which is
homeomorphic to the $2$-dimensional sphere $S^2$.  For this reason,
the volume in this regime grows as $r^2$.

\begin{pro}\label{collapse-lower} There is a constant $c$ such that, uniformly in $a_1\leqslant a_2\leqslant a_3$,
\begin{equation*}
V_{(a_1, a_2, a_3)}(r) \geqslant c a_2^{-2}r^2
  \quad \text{ for }  a_1 \leqslant r\leqslant a_2.
\end{equation*}
\end{pro}

\begin{proof} As usual, it suffices to take $a_2 =1$ (see Remark \ref{scaling-remark}). We proceed along the lines similar to the proof of Proposition \ref{heisenberg-lower}.  For $\eta \in [0,1]$, let $F^\eta : \mathbb{R}^3 \to \SU(2)$ be defined by
  \begin{equation*}
    F^\eta(s_1, s_2, s_3) = \exp(s_1 e_1) \exp(s_2 e_2) H(\eta, s_3),
  \end{equation*}
where $H$ is as in \eqref{H-def}.  Let $J^\eta$ be the Jacobian determinant of $F^\eta$, normalized as in Notation \ref{J-notation}.  Then by the same arguments as in Lemmas \ref{F-lemma} and \ref{J-lemma}, there is a neighborhood $U$ of $(0,0,0) \in \mathbb{R}^3$ and a jointly continuous $j: U \times [0,1] \to \mathbb{R}$ such that
  \begin{equation}\label{J-eta-factor}
    J^\eta(s_1, s_2, s_3) = \eta j(s_1, s_2, s_3, \eta), \qquad (s_1,
    s_2, s_3) \in U, \quad \eta \in [0,1].
  \end{equation}
  We can also directly compute
  \begin{equation*}
    \partial_{s_1} F^\eta(0,0,0) = e_1 \text{ and } \partial_{s_2}
    F^\eta(0,0,0) = e_2.
  \end{equation*}
  For the partial derivative with respect to $s_3$, we can use either
  Lemma \ref{masha-tech} or \eqref{H-basis-formula} to compute
  \begin{equation*}
    \partial_{s_3} F^\eta(0,0,0) = \partial_v H(\eta, 0) = (1 - \cos
    \eta) e_2 + (\sin \eta) e_3.
  \end{equation*}
  Thus, letting $\omega$ be the Riemannian volume form of the bi-invariant metric $g_{(1, 1, 1)}$, we
  have
  \begin{equation*}
    J^\eta(0,0,0) = \frac{1}{16 \pi^2} \omega(e_1, e_2, (1-\cos \eta) e_2 + (\sin \eta)
    e_3) = \frac{1}{16 \pi^2} \sin \eta.
  \end{equation*}
In particular, from \eqref{J-eta-factor}, we have $j(0,0,0,\eta) = \frac{1}{16 \pi^2}\frac{\sin \eta}{\eta} > 0$ for all $\eta \in [0,1]$.  We can thus find a neighborhood $W \subset U$ of $(0,0,0)$ such that $j$ is bounded away from $0$ on $W \times [0,1]$, which implies $J^\eta \geqslant c \eta$ for some constant $c$.

Now choose $\eta>0$ sufficiently small so that $[-\eta, \eta]^3 \subset
W$.  Suppose $\eta a_1 \leqslant r \leqslant \eta$ and set
  \begin{equation*}
    K = \pminterval{\eta} \times \pminterval{r} \times \pminterval{r}.
  \end{equation*}
Note that $K \subset [-\eta, \eta]^3 \subset W$.  Hence, we have $\mu_0(F^\eta(K)) \geqslant c \eta \cdot  m(K) = 8c \eta^2 r^2$.

Also, we have
  \begin{align*}
    d_{(a_1, 1, a_3)}(e, F^\eta(s_1, s_2, s_3)) &\leqslant a_1 s_1 + s_2 + 2
    a_1 \eta + 2 s_3 \\
    &\leqslant a_1 \eta + r + 2 a_1 \eta + 2r \\
    &\leqslant 6r
  \end{align*}
recalling that $r \geqslant a_1 \eta$.  So $F^\eta(K) \subset B_{(a_1, 1, a_3)}(6r)$.  We have thus shown
  \begin{equation*}
    V_{(a_1, 1, a_3)}(6r) \geqslant 8c \eta^2 r^2, \qquad a_1 \eta \leqslant r \leqslant \eta.
  \end{equation*}
  Repeating this argument with $r,\eta$ replaced by $r/6, \eta/6$ (which is valid since we
  still have $[-\eta/6, \eta/6]^3 \subset W$), we have
   \begin{equation*}
    V_{(a_1, 1, a_3)}(r) \geqslant c \eta^2 r^2, \qquad a_1 \eta \leqslant r \leqslant
    \eta
  \end{equation*}
  where a factor of $8/6^4$ has been absorbed into the constant $c$.
  This is the desired result for $a_1 \leqslant r \leqslant \eta$.  For $\eta \leqslant r \leqslant 1$, simply note
  that
  \begin{equation*}
    V_{(a_1, 1, a_3)}(r) \geqslant V_{(a_1, 1, a_3)}(\eta) \geqslant c
    \eta^4 \geqslant c \eta^4 r^2
  \end{equation*}
  and so we have the desired result for all $a_1 \leqslant r \leqslant 1$.
\end{proof}

For the corresponding upper bound, we show that the ball $B_{(a_1, 1, a_3)}(r)$ is contained in a tubular neighborhood of the circle $S = \{\exp(s e_1): s \in \mathbb{R}\}$.

\begin{lem}\label{close-to-circle}
Let $0 < a_1 \leqslant 1 \leqslant a_3 < \infty$.  For any $x \in
\SU(2)$, we may write $x = \exp(s e_1) y$ where $d_{(1,1,1)}(e,y) \leqslant
d_{(a_1, 1, a_3)}(e,x)$.  In particular, $d_{(1,1,1)}(S, x) \leqslant
d_{(a_1, 1, a_3)}(e,x)$.
\end{lem}

\begin{proof}
Note first that without loss of generality we can assume $a_3 = 1$,
since $d_{(a_1, 1, 1)}(e,x) \leqslant d_{(a_1, 1, a_3)}(e,x)$.

Fix $\epsilon > 0$.  Consider the smooth map $\Theta : \mathbb{R}^3 \to \SU(2)$
defined by
\begin{equation*}
  \Theta(z_1, z_2, z_3) = \exp(z_1 e_1 / a_1) \exp(z_2 e_2 + z_3 e_3).
\end{equation*}
Then $d \Theta$ is an isomorphism at $(0,0,0)$, so that $\Theta$ is a
diffeomorphism near $(0,0,0)$.  If we equip $\mathbb{R}^3$ with the
standard Euclidean metric and $\SU(2)$ with the $g_{(a_1, 1,1)}$
metric, then $d\Theta^{-1}_e : T_e \SU(2) \to T_{(0,0,0)}
\mathbb{R}^3$ is an isometry; in particular the operator norm is
$\|d\Theta^{-1}_e\|_{(a_1, 1,1)} = 1$.  Hence we may find some
neighborhood $V$ of $e \in \SU(2)$ such that
$\|d\Theta^{-1}_x\|_{(a_1, 1,1)} \leqslant 1+\epsilon$ for all $x \in V$.
Taking $V$ smaller if necessary, we may also assume that $V$ is a
$g_{(a_1,1,1)}$-normal neighborhood of $e$; that is, for any $x \in V$
there is a $g_{(a_1,1,1)}$-minimizing geodesic from $e$ to $x$
contained in $V$.  Then $\Theta^{-1}$ is a $(1+\epsilon)$-Lipschitz
map from $(V, d_{(a_1, 1,1)})$ into $\mathbb{R}^3$.  So for $x \in V$,
if we write $(z_1, z_2, z_3) = \Theta^{-1}(x)$, we have
\begin{align*}
  x = \Theta(z_1, z_2, z_3) &= \exp(z_1 e_1 / a_1) \exp(z_2 e_2 + z_3
  e_3)
  \\
  &= \exp(s e_1) \exp(t (\cos (\theta) e_2 + \sin(\theta) e_3))
\end{align*}
where we let $s = z_1 / a_1$, $z_2 = t \cos \theta$, $z_3 = t \sin
\theta$.
Moreover,
\begin{align*}
 |t| = \sqrt{z_2^2 + z_3^2} \leqslant |(z_1, z_2, z_3)| \leqslant (1+\epsilon) d_{(a_1, 1, 1)}(e,x).
\end{align*}

Now let $x \in \SU(2)$ be arbitrary.  Let $\gamma : [0,1] \to \SU(2)$
be a $g_{(a_1,1,1)}$-minimizing geodesic from $e$ to $x$, parametrized
by arc length.  For an integer $N$ to be chosen later, let $\tau_i = i/N$ and $x_i =
\gamma(\tau_{i-1})^{-1} \gamma(\tau_i)$, $i=0,\dots,N$, so that $x =
\prod_{i=1}^N x_i$.  Note that by the left invariance of the metric,
\begin{equation*}
  d_{(a_1, 1, 1)}(e, x_i) = d_{(a_1, 1, 1)}(\gamma(\tau_{i-1}),
  \gamma(\tau_i)) = \frac{1}{N} d_{(a_1, 1, 1)}(e,x)
\end{equation*}
since $\gamma$ was parametrized by arc length.  We may now choose $N$
so large that $x_i \in V$ for every $i$.  Then, as above, each $x_i$ may be
written as
\begin{equation*}
  x_i = \exp(s_i e_1) \exp(t_i (\cos (\theta_i) e_2 + \sin(\theta_i) e_3))
\end{equation*}
where
\begin{equation*}
  |t_i| \leqslant (1+\epsilon)d_{(a_1,1,1)}(e, x_i) =
  \frac{1+\epsilon}{N} d_{(a_1,1,1)}(e, x).
\end{equation*}

By repeated application of \eqref{masha-tech-theta}, we may now write
\begin{align*}
  x &= \prod_{i=1}^N \exp(s_i e_1) \exp(t_i (\cos (\theta_i) e_2 +
  \sin(\theta_i) e_3)) \\
  &= \exp(s e_1) \prod_{i=1}^N \exp(t_i (\cos(\phi_i) e_2 +
  \sin(\phi_i) e_3))
\end{align*}
where
\begin{align*}
  s = s_1 + \dots + s_N, \qquad
  \phi_i = \theta_i - s_{i+1} - \dots - s_{N}.
\end{align*}
Setting $y = \prod_{i=1}^N \exp(t_i (\cos(\phi_i) e_2 +
\sin(\phi_i) e_3))$, we have by left-invariance of $d_{(1,1,1)}$ that
\begin{align*}
  d_{(1,1,1)}(e,y) &\leqslant \sum_{i=1}^N d_{(1,1,1)}(e, \exp(t_i (\cos(\phi_i) e_2 +
  \sin(\phi_i) e_3))) \\
  &\leqslant \sum_{i=1}^N |t_i| \\
  &\leqslant (1+\epsilon) d_{(a_1, 1,1)}(e,x).
\end{align*}

To remove the $\epsilon$, we note that for each $n$, we can write $x =
\exp(s_n e_1) y_n$ where, without loss of generality, $s_n \in [-2\pi,
  2\pi]$, and $y_n \in \SU(2)$ with $d_{(1,1,1)}(e,y_n) \leqslant (1 +
\frac{1}{n}) d_{(a_1, 1,1)}(e,x)$.  Since $[-2\pi,2\pi]$ and $\SU(2)$
are compact, we can pass to a subsequence so that $s_n \to s$ and $y_n
\to y$ for some $s,y$, which will then be as desired.
\end{proof}

\begin{pro}\label{collapse-upper}
There is a constant $C$ such that, uniformly in $a_1\leqslant
a_2\leqslant a_3$,
\[
V_{(a_1, a_2, a_3)}(r) \leqslant C a_2^{-2}r^2 \quad \text{ for }  0 \leqslant r\leqslant a_2.
\]
\end{pro}

\begin{proof}
As usual we assume $a_2 = 1$.  Let $K = B_{(1,1,1)}(S, r)$, so that
by the previous lemma $B_{(a_1, 1, a_3)}(r) \subset K$.  It only
remains to estimate the volume of $K$.  Let $N = \lceil
\frac{4\pi}{r} \rceil$, so that $ \frac{4\pi}{r} \leqslant N
\leqslant \frac{4\pi}{r} + 1 \leqslant (4\pi + 1) \frac{1}{r}$.  Set
$x_i = \exp(4 \pi i e_1 / N)$ for $0 \leqslant i \leqslant N$, so that
$x_0 = x_N = e$ and $d_{(1,1,1)}(x_i, x_{i+1}) \leqslant
\frac{4\pi}{N} \leqslant r$.  As
such, the balls $B_{(1,1,1)}(x_i, 2r)$, $1 \leqslant i \leqslant N$,
cover $K$.  Since $(\SU(2), g_{(1,1,1)})$ is a compact 3-dimensional
Riemannian manifold and $\mu_0$ is (up to a constant) its volume
measure, there is a constant $C$ such that $\mu_0(B_{(1,1,1)}(x,R))
\leqslant C R^3$ for any $R$.  So we conclude
  \begin{equation*}
    \mu_0(B_{(a_1, 1, a_3)}(r)) \leqslant \mu_0(K) \leqslant C N (2r)^3 \leqslant 2^3 (4\pi + 1) C r^2.
  \end{equation*}
\end{proof}

\section{Combining the cases}\label{s.Combining}

Combining the foregoing bounds yields the estimates on $V_g(r)$ of Theorem \ref{thV}.

\begin{proof}[Proof of Theorem \ref{thV}]
Similarly to the $V_{(a_1, a_2, a_3)}$ notation, set
\begin{equation}\label{Vbar-a-def}
    \overline{V}_{(a_1, a_2, a_3)}(r) =
    \begin{cases}
      (a_1 a_2 a_3)^{-1} r^3, & 0 \leqslant r \leqslant a_1 a_2/a_3 \\
      (a_1 a_2)^{-2} r^4, & a_1 a_2/a_3 \leqslant r \leqslant a_1 \\
      a_2^{-2} r^2, & a_1 \leqslant r \leqslant a_2 \\
      1, & r \geqslant a_2.
    \end{cases}
  \end{equation}
We need to show that $b_1 \overline{V}_{(a_1, a_2, a_3)}(r) \leqslant V_{(a_1, a_2, a_3)}(r) \leqslant b_2 \overline{V}_{(a_1, a_2, a_3)}(r)$ for some constants $b_1, b_2$ not depending on
 $a_1, a_2, a_3$.  This will establish Theorem \ref{thV} for metrics of the form $g = g_{(a_1, a_2, a_3)}$, recalling from Notation \ref{Va-notation} that $V_g(r)$ differs from $V_{(a_1, a_2,   a_3)}(r)$ by a factor of $(16 \pi^2 a_1 a_2 a_3)^{-1}$.  The  general case follows since, as noted in Corollary \ref{isometry}, every $g \in \mathfrak{L}(\SU(2))$ is isometric to some $g_{(a_1,  a_2, a_3)}$.

The upper and lower bounds in the case $0 \leqslant r \leqslant a_1 a_2/a_3$ are covered by Proposition \ref{euclidean-pro}.

For $a_1 a_2/a_3 \leqslant r \leqslant a_1$, the lower bound is shown by Proposition \ref{heisenberg-lower}.  The upper bound is shown by Proposition \ref{heisenberg-upper} for $a_1 a_2/a_3 \leqslant r \leqslant \eta a_1$, where $\eta$ is a certain small constant, so it remains to handle the case $\eta a_1 \leqslant r \leqslant a_1$.  In this case we can apply Proposition \ref{collapse-upper} to obtain
  \begin{equation*}
    V_{(a_1, a_2, a_3)}(r) \leqslant C a_2^{-2} r^2 \leqslant C \eta^{-2} (a_1
    a_2)^{-2} r^4.
  \end{equation*}
For $a_1 \leqslant r \leqslant a_2$, the desired bounds are given by Propositions \ref{collapse-lower} and \ref{collapse-upper}.

For $r \geqslant a_2$, the lower bound follows simply by noting
  \begin{equation*}
    V_{(a_1, a_2, a_3)}(r) \geqslant V_{(a_1, a_2, a_3)}(a_2)
    \geqslant c
  \end{equation*}
from the bound in Proposition \ref{collapse-lower}.  The upper bound $V_{(a_1, a_2, a_3)}(r) \leqslant 1$ is trivial because $V_{(a_1, a_2, a_3)}(r)$ is the volume with respect to the probability measure $\mu_0$.
\end{proof}

To prove Theorem \ref{thM}, it now suffices to show that the function $\overline{V_{g}}$, or equivalently $\overline{V}_{(a_1, a_2, a_3)}$ as in \eqref{Vbar-a-def}, satisfies a uniform volume doubling condition.  This is an elementary calculation which we insert here for convenience.

\begin{lem}\label{Vbar-doubling}
  For any $a_1 \leqslant a_2 \leqslant a_3$, any $r \geqslant 0$ and $k \geqslant 1$, we have
  \begin{equation*}
    \overline{V}_{(a_1, a_2, a_3)}(kr) \leqslant k^4 \, \overline{V}_{(a_1, a_2, a_3)}(r).
  \end{equation*}
\end{lem}

\begin{proof}
 We have ten cases depending on which of the four regions defined in \eqref{Vbar-a-def} are occupied by $r$ and $kr$.

 If $r,kr$ occupy the same region, then the result is immediate.  For instance, when $0 \leqslant r \leqslant kr \leqslant a_1 a_2/a_3$, then we have $\overline{V}(kr) / \overline{V}(r) = k^3$ (we suppress the subscripts).  In the other similar cases, we get $k^4$, $k^2$ or $1$; all are bounded by $k^4$.

The next cases are when they occupy consecutive regions.
  \begin{itemize}
  \item If $0 \leqslant r \leqslant a_1 a_2/a_3 \leqslant kr \leqslant
    a_1$, then $\overline{V}(kr) / \overline{V}(r) = k^4 (a_1 a_2/a_3)^{-1} r \leqslant k^4$ because
    $r \leqslant a_1 a_2/a_3$.
  \item If $a_1 a_2 / a_3 \leqslant r \leqslant a_1 \leqslant kr
    \leqslant a_2$, then $\overline{V}(kr) / \overline{V}(r) = k^2
    a_1^2 r^{-2} \leqslant k^4$, using $a_1 \leqslant kr$.
  \item If $a_1 \leqslant r \leqslant a_2 \leqslant kr$, then
    $\overline{V}(kr) / \overline{V}(r) = a_2^2 r^{-2} \leqslant k^2
    \leqslant k^4$, using $a_2 \leqslant kr$.
  \end{itemize}

  The remaining cases follow by combining those already shown.  For
  instance, if $0 \leqslant r \leqslant a_1 a_2/a_3 \leqslant
  a_1 \leqslant kr \leqslant a_2$, choose $1 \leqslant k' \leqslant k$ so that
  $a_1 a_2/a_3 \leqslant k'r \leqslant a_1$.  Then by the previous
  cases we have
  \begin{equation*}
    \overline{V}(kr) \leqslant \left(\frac{k}{k'}\right)^4
    \overline{V}(k'r) \leqslant k^4 \overline{V}(r).
  \end{equation*}
  The last two cases are similar.
\end{proof}

Combining Theorem \ref{thV} and Lemma \ref{Vbar-doubling} (with $k=2$) establishes
Theorem \ref{thM}, with $D = 16 b_2/b_1$.

\section{Diameter bounds}\label{s.Diameter}

In this brief section, we prove the remark following Theorem \ref{thV}: for any metric $g \in
\mathfrak{L}(\SU(2))$, the diameter $\operatorname{diam}_g(\SU(2))$ is uniformly comparable to $a_2$, the square root of the middle eigenvalue.

An interesting consequence is that, by inspection of \eqref{ricci-bound}, there is no uniform lower bound on the Ricci curvatures of the metrics $g \in \mathfrak{L}(\SU(2))$, even after
rescaling to constant diameter; the metrics $g_{(1,1,a_3)}$, as $a_3 \to \infty$, have comparable diameters, but their Ricci curvatures in the $e_3$ direction tend to $-\infty$.  As such, the uniform volume doubling bound of Theorem \ref{thM} cannot be obtained solely by Ricci curvature considerations as in Section \ref{s.EuclideanRegime}.

\begin{pro}\label{diameter-bounds}
  For a left-invariant Riemannian metric $g \in \mathfrak{L}(\SU(2))$,
  let $a_2$ be the square root of the middle eigenvalue of the matrix
  $A_g$, as in Theorem \ref{thV}.  There are universal constants
  $0<D _0 \leqslant D_\infty<+\infty$ such that
\[
D_0 a_2\leqslant \operatorname{diam}_g(\SU(2))\leqslant D_\infty a_2.
\]
\end{pro}

\begin{proof} By Corollary \ref{isometry}, we can assume without loss
  of generality that $g = g_{(a_1, a_2, a_3)}$ for some $a_1 \leqslant
  a_2 \leqslant a_3$, and by scaling, we can assume $a_2=1$.

  For an upper bound, we consider a sub-Riemannian metric on $\SU(2)$.
  Let $\mathcal{H} \subset T \SU(2)$ be the two-dimensional sub-bundle
  spanned at each point by the left translates of $\hat{e}_1,
  \hat{e}_2$, and let $g_{( 1, 1,\infty)}$ be the left-invariant
  sub-Riemannian metric on $\mathcal{H}$ making $\hat{e}_1, \hat{e}_2$
  orthonormal.  Then $(\SU(2), \mathcal{H}, g_{(1,1,\infty)})$ is a
  sub-Riemannian manifold.  The sub-bundle $\mathcal{H}$ satisfies
  H\"ormander's bracket-generating condition, since $[\hat{e}_1,
    \hat{e}_2] = \hat{e}_3$, and so by the Chow--Rashevskii theorem
  \cite[p. 43]{MontgomeryBook2002}, the sub-Riemannian (or
  Carnot--Carath\'eodory) distance $d_{(1, 1, \infty)}$ is finite and
  induces the original manifold topology.  Since $\SU(2)$ is compact,
  it has finite diameter under $d_{(1,1,\infty)}$.  Let $D_\infty$ be
  this diameter.  It is clear that for any $v \in T\SU(2)$, we have
  $g_{(a_1, 1, a_3)}(v,v) \leqslant g_{(1, 1, \infty)}(v,v)$ (where
  for $v \notin \mathcal{H}$ we can take $g_{(1, 1, \infty)}(v,v) =
  \infty$), so the same inequality holds for their distances, and we
  have shown that the diameter under $g_{(a_1, 1, a_3)}$ is bounded
  above by $D_\infty$.

For the lower bound, consider the pseudo-metric $g_{(0,1,1)}$ for
which $g_{(0,1,1)}(\hat{e}_1, \hat{e}_1) = 0$ and $\hat{e}_2,
\hat{e}_3$ are orthonormal.  Then the pseudo-distance $d_{(0,1,1)}$ is
symmetric and satisfies the triangle inequality, but is not positive
definite.  For instance, $d_{(0,1,1)}(e, \exp(s \hat{e}_1)) = 0$ for
any $s$.  However, we claim $d_{(0,1,1)}$ is not identically zero, so
that $\SU(2)$ has nonzero diameter under $d_{(0,1,1)}$.  As above,
$d_{(0,1,1)}$ is a lower bound for any $d_{(a_1, 1, a_3)}$, so we may
take $D_0$ to be the $d_{(0,1,1)}$-diameter of $\SU(2)$.

Indeed, let $S = \{\exp(s \hat{e}_1) : s \in \mathbb{R}\}$ be the
subgroup generated by $\hat{e}_1$.  Suppose $d_{(0,1,1)}(e,x) = 0$; we
claim that $x \in S$.  For any $\epsilon > 0$, we can choose $a$ so
small that $d_{(a,1,1)}(e,x) < \epsilon$.  By Lemma
\ref{close-to-circle}, we can write $x = \exp(s e_1) y$ where
$d_{(1,1,1)}(e,y) < \epsilon$.  Thus $d_{(1,1,1)}(S, x) <
\epsilon$.  Since $\epsilon$ was arbitrary and $S$ is closed, we
conclude that $x \in S$.  So for any $x \in \SU(2) \setminus S$, we
have $d_{(0,1,1)}(e,x) > 0$.
\end{proof}

\begin{rem}
  In effect, the pseudo-metric space $(\SU(2), d_{(0,1,1)})$ is the
  two-di\-men\-sion\-al left coset space $\SU(2) / S$, which is
  homeomorphic to $S^2$.  This statement is not so obvious as it might
  appear.  For instance, suppose we instead consider the Heisenberg
  group $\mathbb{H}^3$ with the standard basis $\{X,Y,Z\}$ for
  $\mathfrak{h}^3$ satisfying $[X,Y]=Z$, $[X,Z]=[Y,Z]=0$, and a
  left-invariant pseudo-metric $g$ with $g(X,X)=0$ and $Y,Z$
  orthonormal.  Then the resulting pseudo-metric space is only
  one-dimensional, and in particular it does not equal the quotient of
  $\mathbb{H}^3$ by $\{\exp(tX):t \in \mathbb{R}\}$.  Indeed, by
  writing $\exp(s^2 Z) = \exp(s\epsilon^{-1} X) \exp(s \epsilon Y)
  \exp(-s \epsilon^{-1} X) \exp(-s \epsilon Y)$ where $\epsilon \to
  0$, we see that we can reach the $z$-axis by paths of arbitrarily
  small length with respect to this metric, by making a rectangle that
  is very large in the $X$ direction and very small in $Y$.  However,
  compactness prevents this phenomenon in $\SU(2)$.
\end{rem}

\begin{rem}\label{diam-explicit}
In a recent article \cite{Podobryaev2018}, A.~V.~Podobryaev has  computed the diameter of the metrics $g_{(a_1, a_2, a_3)}$ in the case where two of the three parameters $a_1, a_2, a_3$ are equal.   This leads to the explicit values $D_0 = \pi$, $D_\infty = 2\pi$.  The value $D_\infty = 2\pi$ also follows from a sub-Riemannian  distance formula proved in \cite{BoscainRossi2008a}.
\end{rem}

\section{Consequences of volume doubling}\label{s.Consequences}

Let $\left( M, g \right)$ be a Riemannian manifold, and $\Delta_{g}$
the (positive) Laplace--Beltrami operator associated with the metric $g$. The gradient $\nabla_{g}$ is determined by the metric $g$ and we let
 \[
 \vert \nabla_{g} f\vert_{g}^{2}:=g\left( \nabla_{g} f, \nabla_{g} f \right).
 \]

The connection between the Laplace--Beltrami operator and the gradient
is given by
 \[
 \int_{M} f \Delta_{g} f d\mu_{g}=\int_{M} \vert \nabla_{g} f\vert^{2}d\mu_{g},
 \]
where as before $\mu_{g}$ is the Riemannian volume measure. Finally the heat kernel is the fundamental solution to the heat equation with the Laplace--Beltrami operator $\Delta_{g}$, which equivalently can be described as the kernel for the heat semigroup
\[
P_{t}f\left( x \right)=e^{-t\Delta_{g}}f\left( x \right)= \int_{M} f\left( y \right)p_{t}^{g}\left( x, y \right)d\mu_{g}\left( y \right).
\]

We concentrate on the case when $M$ is a compact Lie group. Namely, let $K$ be a connected compact group equipped with a left-invariant Riemannian metric $g \in \mathfrak{L}(K)$. In this case, the heat kernel  $p_{t}^{g}\left( x, y \right)$ is a symmetric function of $(x,y)$ and is invariant under left multiplication, that is, $p_{t}^{g}\left( x, y \right) = p^g_t(e,x^{-1}y)=p^g_t(e,y^{-1}x)$.  Abusing notation, we write $p^g_t(z):=p^g_t(e,z)$. In addition, the heat kernel satisfies the Chapman--Kolmogorov equations
\begin{equation}\label{e.8.1}
p_{s+t}^{g}\left( x \right)=\int_{K} p_{s}^{g}\left(y^{-1}x\right) p_{t}^{g}\left( y \right) d\mu_{g}, \hskip0.1in s, t >0,
\end{equation}
which implies (using symmetry and multiplication invariance) that
\begin{equation}\label{e.8.2}
p_{s+t}^{g}\left( e \right)=\int_{K} p_{s}^{g}\left(y\right) p_{t}^{g}\left( y \right) d\mu_{g}.
\end{equation}
As mentioned before, the volume doubling constant  is quantitatively related to many analytic properties of the Laplace--Beltrami operator $\Delta_g$. Given a Riemannian metric $g$ on a compact manifold $M$, let
\begin{equation}\label{e.spectrum}
0=\lambda_{g, 0} < \lambda_{g}= \lambda_{g, 1} \leqslant \cdots\leqslant \lambda_{g, i}\leqslant \cdots
\end{equation}
be the eigenvalues of $\Delta_g$, repeated according to multiplicity. In the case when $M$ is a compact Lie group, we use will use repeatedly the following connection between the heat kernel $p_{t}^{g}$ and the eigenvalues
\begin{equation}\label{e.HeatKernelEigenvalue}
\mu_{g}(K) p^g_t(e)=\mu_{g}(K) p^g_t(x, x)= \sum _{i=o}^{\infty} e^{-t\lambda_{g,i}}.
\end{equation}
We discuss some of the properties of $\Delta_g$ and  the volume doubling constant here, including a spectral gap, Weyl eigenvalue counting function, parabolic Harnack inequalities, and heat kernel bounds.

\begin{defin}\label{d.UniformlyDoubling} Let $K$ be a connected real Lie group. We say that $K$ is \textbf{uniformly doubling with constant at most $D$} if there is a constant $D$ such that
\[
D_{g} \leqslant D
\]
for all left-invariant metrics Riemannian metrics $g \in \mathfrak{L}(K)$.
\end{defin}
Observe that by Lemma \ref{Vbar-doubling}  we see that on $\SU(2)$
\[
\frac{V_{g}\left( r \right)}{V_{g}\left( s \right)}\leqslant D \left(\frac{r}{s}\right)^{4}  \text{ for any } 0 < s \leqslant r.
\]
This can be compared with a more general statement as follows. Suppose $\left( X, d, \mu \right)$ is a metric measure space, then one can ask if there are constants $D^{\prime}>0$ and $\delta>0$ such that for any $0 < s \leqslant r$ and $x \in X$
\begin{equation}\label{e.DoublingComparison}
\frac{V\left( x, r \right)}{V\left( x, s \right)}\leqslant D^{\prime} \left(\frac{r}{s}\right)^\delta,
\end{equation}
where $V\left( x, r \right):=\mu\left( B\left( x, r \right)\right)$. If the metric measure space $\left( X, d, \mu \right)$ is doubling with constant at most $D$, then by \cite[Section 4.2]{CoulhonSikora2008} and \cite[Lemma 5.2.4]{Saloff-CosteBook2002}  we see that \eqref{e.DoublingComparison} holds with $D^{\prime}=D$ and $\delta=\frac{\ln D}{\ln 2}$.
Indeed, if $\lfloor \cdot \rfloor$ denotes the integer part of a real number (floor function), we see that
\begin{align*}
& V\left( x, r \right) \leqslant V\left( x, 2^{\lfloor {\frac{\ln \left( \frac{r}{s}\right)}{\ln 2}} \rfloor +1}s \right) \leqslant D^{1+\frac{\ln \left( \frac{r}{s}\right)}{\ln 2}}V\left( x, s \right) = D \left( \frac{r}{s}\right)^{\frac{\ln D}{\ln 2}}V\left( x, s \right).
\end{align*}

Below we state several interesting properties which would follow from Conjecture \ref{MC1}.  First and foremost, we note that it implies a uniform version of the Poincar\'e inequality for metric balls stated in Corollary \ref{uniform-poincare}.  This is the key to a host of other consequences. In particular, by Theorem \ref{thM} these properties hold on $\SU(2)$. In some instances, Theorem \ref{thV} provides a particularly explicit form of these statements.

\subsection{The Poincar\'e inequality on compact Lie groups}\label{a.Poincare}

The following theorem is proved in \cite[Section 5.6.1]{Saloff-CosteBook2002}.  The first instance of this type of inequality appeared in \cite{Varopoulos1987a}; a discrete version of this inequality is one of the key elements of B.~Kleiner's proof of Gromov's theorem on groups of polynomial growth \cite{Kleiner2010}.

\begin{theo}\label{poincare-thm}
  Let $K$ be a compact Lie group equipped with a left-invariant   Riemannian metric $g$.  On any ball $B_g(x,r)$, we have the   Poincar\'e inequality
\begin{equation}\label{poincare-eqn}
 \int_{B_{g}\left(x,r \right)}|f-f_{x, r}|^2 d\mu_{g} \leqslant 2
 r^{2} D_{g} \int_{B_{g}\left(x,2r \right)} |\nabla_{g} f|_{g}^2
 d\mu_{g} \text{ for all } f \in \mathcal C^\infty(B_{g}\left(x,2r
 \right)),
 \end{equation}
 where $f_{x, r}:=\int_{B_{g}\left(x,r \right)}f d\mu_{g}$ denotes the mean of $f$ over $B_{g}\left(x,r \right)$, and $D_g$ is the volume
 doubling constant of $(K,g)$.
\end{theo}

\begin{cor}\label{uniform-poincare}
If $K$ is uniformly  volume doubling with constant at most $D$, then the Poincar\'e inequality \eqref{poincare-eqn} holds with the same constant $D$ for  every $g \in \mathfrak{L}(K)$.  In particular, by Theorem \ref{thM} this is true for $K = \SU(2)$.
\end{cor}

The proof given in \cite{Saloff-CosteBook2002} also establishes, by a
straightforward modification, the following $L^p$ Poincar\'e
inequality:
\begin{theo}\label{Lp-poincare}
  In the same notation as Theorem \ref{poincare-thm}, for any $1
  \leqslant p < \infty$, we have
  \begin{equation}\label{Lp-poincare-eqn}
  \int_{B_{g}\left(x,r \right)}|f-f_{x, r}|^p d\mu_{g} \leqslant
  (2r)^{p}
  D_{g} \int_{B_{g}\left(x,2r \right)} |\nabla_{g} f|_{g}^p
 d\mu_{g} \text{ for all } f \in \mathcal C^\infty(B_{g}\left(x,2r
 \right)).
  \end{equation}
\end{theo}
(In the special case $p=2$, one can improve the constant by a factor
of $\frac{1}{2}$ to recover \eqref{poincare-eqn}.)

Note that the weak Poincar\'e inequality \eqref{poincare-eqn}  and volume doubling imply that the strong Poincar\'e inequality holds, that is, \eqref{poincare-eqn} with the same ball $B_{g}\left(x,r
\right)$ on both sides (and the same for the $L^p$ Poincar\'e
inequality \eqref{Lp-poincare-eqn}). This is shown by a covering argument; see \cite{Jerison1986a} and \cite[Section 5.3.2]{Saloff-CosteBook2002}.  In particular, this implies that on  a uniformly  doubling group $K$ the lowest eigenvalue  $\lambda_{N,g,r}$ of the Laplacian $\Delta_g$ with Neumann boundary condition on the ball $B_{g}(x,r)$ satisfies $c r^{-2}\leqslant \lambda_{N,g,r}\leqslant C r^{-2}$ uniformly over all $g \in \mathfrak{L}(K)$ and $r\in (0,\operatorname{diam}_g]$.

\subsection{Spectral gap}\label{ss.spectral-gap}
Let $\lambda_{g}$ be the lowest non-zero eigenvalue for the
Laplace-Beltrami operator $\Delta_{g}$.  We show that when $K$ is
uniformly doubling, we obtain the uniform upper bound \eqref{e.1.0}
for $\lambda_g$, matching the lower bound \eqref{li-eigenvalue-lower}
obtained in \cite{LiP1980a}, up to a constant depending on
the doubling constant.

\begin{theo}\label{spectral-gap} Assume that the compact  connected
  Lie group $K$ is uniformly doubling with constant at most $D$.
  For any metric $g\in \mathfrak{L}(K)$, the lowest non-zero
  eigenvalue $\lambda_{g}$ of the Laplacian $\Delta_g$ satisfies
\[
\frac{\pi^2}{4 \operatorname{diam}_g^2}\leqslant  \lambda_{g}
\leqslant  \frac{16 D^2}{\operatorname{diam}_g^2}.
\]
 \end{theo}

\begin{proof}
  As mentioned earlier, the lower bound was proved in
  \cite{LiP1980a}. (An improved lower bound was recently obtained in
  \cite{JudgeLyons2017}.) To obtain an upper bound, we note that
 \begin{equation}\label{e.VariationalLambda}
 \lambda_{g}= \min\left\{\frac{\int |\nabla_g f|^2 d\mu_g}{\|f\|_2^2}: f\neq 0, \hskip0.05in   \int_K fd\mu_{g}=0, \hskip0.05in  f\in \operatorname{Lip}(K) \right\},
 \end{equation}
 We construct an
  appropriate test function to use in \eqref{e.VariationalLambda}. Let
  $y$ be a point which realizes the diameter of $K$ under $g$, i.e.,
  $d_g(e,y)=\operatorname{diam}_g$. Let $R=\operatorname{diam}_g/2$.
  For any $z\in K$, let $f_{z,r}(x)=(r-d_g(z,x))_+$ be the tent
  function over the ball $B_g(z,r)$; observe that this is a Lipschitz
  function with gradient $|\nabla_g f_{z,r}|\leqslant 1$ (almost
  everywhere).  As a test function, take $f_R= f_{e,R}-f_{y,R}$. By
  group invariance, $\int_K f_Rd\mu_g=0$ and $\int_K |\nabla_g
  f_R|^2d\mu_g \leqslant \mu_g(K)=V_g(2R)$.  To estimate the
  $L^2$-norm of $f_R$ from below, observe that $|f_R|$ is at least $R/2$
  on two disjoint balls of radius $R/2$. Hence we have
  $\|f_R\|_2^2\geqslant (R/2)^2V_g(R/2)$. Plugging this in the
  variational formula \eqref{e.VariationalLambda} for $\lambda_{g}$
  yields
\[
\lambda_{g} \leqslant  \frac{4V_g(2R)}{R^2V_g(R/2)}\leqslant \frac{16 D^2}{\operatorname{diam}_g^2}.
\]
\end{proof}

In the special case when $K = \SU(2)$, we have from Proposition \ref{diameter-bounds} that the diameter $\operatorname{diam}_g(\SU(2))$ is uniformly comparable to the parameter $a_2$ of $g$ (as defined in Notation \ref{n.MetricsParameters}), and hence we have the following
statement.
\begin{cor} There are positive constants $0< c \leqslant C<\infty$ such that for all $g\in \mathfrak{L}(\SU(2))$ with parameters $0<a_1\leqslant a_2 \leqslant a_3$ as in Notation \ref{n.MetricsParameters}, we have
\begin{equation}\label{spectral-gap-SU2}
  \frac{c}{a_2^2}\leqslant  \lambda_{g}\leqslant  \frac{C}{a_2^2}.
\end{equation}
\end{cor}

\begin{rem}
  In a very recent preprint \cite{Lauret2018}, E.~A.~Lauret has given an
  exact expression for the smallest eigenvalue $\lambda_g$ of
  $\SU(2)$ in terms of the parameters of the metric, which in our
  notation reads as follows:
  \begin{equation}\label{spectral-gap-lauret}
    \lambda_g = \min\left\{\frac{1}{4}\left(\frac{1}{a_1^{2}}  +
    \frac{1}{a_2^{2}} + \frac{1}{a_3^{2}}\right), \frac{1}{a_2^{2}} + \frac{1}{a_3^{2}}\right\}.
  \end{equation}
  Indeed, \eqref{spectral-gap-lauret} is consistent with
  \eqref{spectral-gap-SU2}.

  In earlier work, as part of a more general construction, H.~Urakawa \cite{Urakawa1979} computed $\lambda_g$ for a particular one-parameter family of metrics $g(t)$ on $\SU(2)$, which in our
  notation is
  \begin{equation*}
    g(t) =
    \begin{cases}
      g_{(t/\sqrt{2},\; 1/\sqrt{2t},\; 1/\sqrt{2t})}, & 0 < t \leqslant 1 \\
      g_{(1/\sqrt{2t},\; 1/\sqrt{2t},\; t/\sqrt{2})}, & 1 \leqslant t
      < \infty.
    \end{cases}
  \end{equation*}
  See \cite[Theorem 5]{Urakawa1979}.  This family has the property  that the volume $\mu_{g(t)}(\SU(2))$ is the same for all $t$, while   $\lambda_{g(t)} \sim t$.

Urakawa's example answered, in the  negative, a previous question of M.~Berger \cite{Berger1973a}: whether   we have $\lambda_g \leqslant C(M) \mu_g(M)^{-2/n}$ on any $n$-dimensional compact connected manifold $M$, with a constant $C(M)$ depending on $M$ but not on the metric $g$.  It is interesting to compare this with Theorem   \ref{spectral-gap}, which implies that, when $M$ is a uniformly   doubling group $K$ and the metrics are left-invariant, the quantity  $\mu_g(K)^{1/n}$ in Berger's statement ought to be replaced with $\operatorname{diam}_g$.
\end{rem}

\subsection{Heat kernel estimates}\label{ss.HeatKernelEst}
In the section we would like to comment on the heat kernel estimates \eqref{e.1.1} for uniformly doubling compact Lie groups. Given a complete Riemannian manifold that satisfies the volume doubling property and the Poincar\'e inequality \eqref{poincare-eqn}, there are several ways to obtain heat kernel upper bounds. One of the most direct and efficient is based on the notion of a Faber--Krahn inequality as  developed in \cite{Grigoryan1994a, Carron1996a} or the equivalent notion of local Sobolev inequality (see \cite[Section 5.2]{Saloff-CosteBook2002}).

Assuming that doubling and the Poincar\'e inequality hold, these methods provide the heat kernel upper bound in terms of the volume
\begin{equation}\label{e.HeatKernelUpper}
p_t(x,y)\leqslant  \frac{C_1(\varepsilon)}{
\sqrt{V(x,\sqrt{t})V(y,\sqrt{ t})}}\exp\left(- \frac{d(x,y)^2}{4(1+\varepsilon)t}\right)
\end{equation}
with a constant $C_1(\varepsilon)$  \cite[Equation(5.2.17)]{Saloff-CosteBook2002} and \cite{Saloff-Coste1992a}  that depends only on $\varepsilon\in (0,1)$ and the constants involved in the doubling property and the Poincar\'e inequality. Here  $V\left( x, r \right)$ denotes the volume of the ball of radius $r>0$ around the point $x$.

In fact, these arguments provide the more precise bound of the type
\[
p_t(x,y)\leqslant  \frac{C_1(1+d(x,y)^2/4t)^{\kappa}}{\sqrt{V(x,\sqrt{t})V(y,\sqrt{t})}} \exp\left(- \frac{d(x,y)^2}{4t}\right)
\]
for some $\kappa>0$. The best value of $\kappa$ that can be obtained from these arguments is $\kappa=\delta/2$,  where $\delta$ is  as in \eqref{e.DoublingComparison}, e.g.  \cite[Section 5.2.3]{Saloff-CosteBook2002},  \cite[Corollary 4.2]{Sturm1996a} and variations on the arguments in \cite{CoulhonSikora2008}. Further,  one also obtains the time derivative estimates such as in \cite[Corollary 2.7]{Sturm1995a})
\[
|\partial_t^kp_t(x,y)|\leqslant  \frac{C_k(1+d(x,y)^2/4t)^{k+\delta/2}}{t^{k}\sqrt{V(x,\sqrt{t})V(y,\sqrt{t})}} \exp\left(- \frac{d(x,y)^2}{4t}\right).
\]
In addition, \cite{CoulhonSikora2008} provides assorted estimates for the heat kernel in complex time and pointers to further references.

The proofs of these estimates simplify, and a greater varieties of arguments can be employed, when the volume of balls is independent of the center, which is the case for left-invariant metrics on Lie groups.

\begin{theo}\label{t.TimeDervativeHK} Let $K$ be a compact Lie group. If $K$ is uniformly volume doubling with constant at most $D$, then for each integer $k=0,1, \dots$ there exists a constant $C_k$ depending only on $D$ and $k$ such that for any $g\in \mathfrak{L}(K)$, and for all $x, y\in K$ and $t>0$
\[
|\partial_t^kp^g_t(x,y)|\leqslant \frac{C_k(1+d_g(x,y)^2/4t)^{k+\delta/2}}{t^kV_g(\sqrt{t})} \exp \left(- \frac{d_g(x,y)^2}{4t}\right),
 \]
 where $\delta$ is as in \eqref{e.DoublingComparison}.
\end{theo}
Regarding a lower bound, the only directly applicable results are proved by a simple  chaining argument using the parabolic Harnack inequality discussed in Section \ref{s.HarnackIneq}.  Assuming that doubling and the Poincar\'e inequality hold, this line of reasoning  provides the following  heat kernel lower bound
\[
p^{g}_t(x,y)\geqslant  \frac{c_2}{V(x,\sqrt{t})} \exp\left(- C_2\frac{d(x,y)^2}{t}\right),
\]
where $0<c_2,C_2$ depends only on the constants involved in the doubling property and the Poincar\'e inequality.  See, for instance, \cite[Section 5.4.6]{Saloff-CosteBook2002}  and \cite[Corollary 4.10]{Sturm1996a}.

\begin{theo}\label{t.HeatKernelLower} Let $K$ be a compact Lie group. If $K$ is uniformly volume doubling with constant at most $D$ then there exist  positive constants $c$ and $A$ depending only on $D$ such that, for any $g\in \mathfrak{L}(K)$,  for all $x, y \in K$ and $t>0$
\[
p^g_t(x,y)\geqslant \frac{c}{V_g(\sqrt{t})} \exp \left(- A\frac{d_g(x,y)^2}{t}\right).
\]
\end{theo}
For $\SU(2)$, Lemma \ref{Vbar-doubling} shows that  we can take  $\delta =4$ (uniformly over $\mathcal L(\SU(2))$) in Theorem \ref{t.TimeDervativeHK}, and this gives the following result.

\begin{theo} There exist constants $0<c, A$ and for each $k=0,1,\dots$, a constant $C_k$,
such that, for all $g\in \mathcal L(\SU(2))$ and all $x, y\in \SU(2)$, $t>0$, we have
\[|\partial_t^kp^g_t(x,y)|\leqslant \frac{C_k(1+d_g(x,y)^2/4t)^{k+2}}{t^kV_g(\sqrt{t})} \exp \left(- \frac{d_g(x,y)^2}{4t}\right)
\]
and
\[
p^g_t(x,y)\geqslant \frac{c}{V_g(\sqrt{t})}\exp \left(- A\frac{d_g(x,y)^2}{t}\right).
\]
\end{theo}
\begin{rem} The results in  \cite{Varopoulos1989a} imply that for each metric $g\in \mathfrak{L}(K)$ (in particular for $\SU(2)$) and $\epsilon\in (0,1)$ there is a constant $c_{\varepsilon \left( g \right)}>0$ such that, for all $x, y, t$
\[
p^g_t(x,y)\geqslant \frac{c_{\varepsilon \left( g \right)}}{V_{g}(\sqrt{t})}\exp \left(- \frac{d_g(x,y)^2}{4(1-\varepsilon)t}\right).
\]
However, it is not clear that the arguments in \cite{Varopoulos1989a} are sufficient to yield a constant $c_\varepsilon$ that is uniform in $g$, even if one assumes that the group $K$ is uniformly doubling. This remains an open question, although we conjecture that this inequality holds uniformly.
\end{rem}
\begin{rem} Detailed asymptotics originally developed by S.~Molchanov in \cite{Molchanov1975} show that for the heat kernel on the $n$-sphere equipped with its canonical round metric and with $x$ and $y$ being antipodal points (e.g., the south and north poles)
\[
p_t(x,y)\sim c_n t^{-n/2}  \left( \frac{d(x,y)^2}{t} \right)^{(n-1)/2}\exp(-d(x,y)^2/4t)
\]
as $t$ tends to $0$. This shows that one cannot dispense entirely with the factor $(1+d_g(x,y)^2/t)^\kappa$ in heat kernel upper bounds, even  on $\SU(2)$. For more on this, see \cite{Neel2007}.
\end{rem}

\subsection{Harnack inequality}\label{s.HarnackIneq} Let the parabolic Harnack constant $H(M, g)$ be the infimum of  all real $H$ such that for any $x \in M$, $r>0$ and any positive solution $u$ of the heat equation on $(M, g)$ in $(s, s+4r^2)\times B(x, 2r)$, it holds that
\begin{equation}\label{e.2.4}
\sup_{Q_-}\{ u\}\leqslant H \inf_{Q_+}\{u\},
\end{equation}
where $Q_-=(s+r^2,s+2r^2)\times B(x,r)$ and
$Q_+=(s+3r^2,s+4r^2)\times B(x,r)$.  In particular, for any
connected compact real Lie group $K$ equipped with $g \in \mathfrak{L}(K)$, we denote by $H(K, g)$ the best constant in the parabolic
Harnack inequality \eqref{e.2.4}. Then one can ask if the parabolic
Harnack inequality is satisfied uniformly over all $g \in \mathfrak{L}(K)$.

\begin{pro}[See \cite{Grigoryan1991, Saloff-Coste1992b}]\label{p.2.3} Let $K$ be a connected compact Lie group. Assume that $K$ is uniformly doubling with constant at most $D$; then there is a $H\left( D \right)$ such that
\[
H(K, g) \leqslant H\left( D \right)
\]
for all $g \in \mathfrak{L}(K)$.
\end{pro}

In particular, Theorem \ref{thM} implies the following.

\begin{cor}[Uniform Harnack inequality for $\SU(2)$]
The parabolic Harnack inequality is satisfied uniformly over all  $g\in \mathfrak{L}(\SU(2))$.
\end{cor}

\subsection{Gradient inequalities}\label{s.GradEst}

In addition to the Harnack inequality \eqref{e.2.4} several related useful inequalities involve gradient estimates. For instance, one can consider the property that for any $x\in M$, $r>0$ and any positive solution $u$ of the heat equation on $(M, g)$ in $(s,s+4r^2)\times B(x,2r)$, it holds that

\begin{equation}\label{e.2.3}
\sup_{Q_-}\{|\nabla_{g} u|_{g} \}\leqslant H_1 t^{-1/2}\inf_{Q_+}\{u\}
\end{equation}
with $Q_{-},Q_{+}$ defined as above. Or one may prefer the Li-Yau parabolic inequality for global positive solutions $u\left( t, x \right)$ of the heat equation on  $(M, g)$ in $(0,T)\times M$,
\begin{equation}\label{LiYau}
|\nabla_{g} \log u|_{g}^2 -\partial_t \log u \leqslant  \frac{H'_1}{t}.
\end{equation}

In this direction, we can only prove the following weaker result for the heat kernel $p_t^g \left( x \right)$.

\begin{theo}\label{t.GradientEstimates} Assume that $K$ is uniformly
  doubling with constant at most $D$. Then there is a constant
  $C\left( D \right)$ such that
\begin{align}
& \vert \nabla_g p_t^g \left( x \right)\vert_{g} \leqslant
\frac{C\left( D \right)}{\sqrt{t} V_{g}\left( \sqrt{t} \right)}\left( 1+\frac{d_g^2(e, x)}{4t} \right)^{3\delta +1}
\exp \left( - \frac{d_g^2(e, x)}{4t} \right),
\label{e.2.5}
\\
&  \Vert \nabla_gp_t^g \Vert_{L^{1}}=\int_K \vert \nabla_g p_t^g \left( x \right)\vert_{g} d\mu_g(x)\leqslant C\left( D \right) t^{-1/2}
\label{e.2.6}
\end{align}
where $\delta=\delta \left( D \right)$ is as in \eqref{e.DoublingComparison}.
\end{theo}

\begin{proof}
Spectral theory easily gives
\begin{equation}\label{e.8.10}
\Vert \Delta_{g} P^g_t \Vert_{L^{2} \rightarrow L^{2}}=\|\partial _t P^g_t\|_{L^{2} \rightarrow L^{2}}=\sup_{\lambda>0}\left\{ \lambda e^{-t\lambda}\right\} \leqslant   \left(et\right)^{-1}\leqslant t^{-1},
\end{equation}
where we used the operator norm on $L^2(K,\mu_g)$ for $P^g_t$. Now observe that by \eqref{e.8.1}
\[
\nabla_g p^g_t(x)= \int_K \nabla_g p^g_{t/2}(y^{-1} x)p^g_{t/2}(y)d\mu_g(y).
\]
Hence
\begin{eqnarray*}
\vert \nabla_g p_t^g(x)\vert_{g}^{2} & \leqslant &  \left(\int_K \vert \nabla _g p^g_{t/2}(y^{-1} x)\vert_{g}p^g_{t/2}(y)d\mu_g(y) \right)^2\\
&\leqslant & \int_K \vert \nabla _g p^g_{t/2}(y)\vert_{g}^{2} d\mu_g(y) \int_K |p^g_{t/2}(y)|^2d\mu_g(y) .
\end{eqnarray*}

By \eqref{e.8.2} we have $\|p^g_{t/2}\|_2^{2}=\int_K |p^g_{t/2}(y)|^2d\mu_g(y)=p^g_t(e)$  and
\begin{eqnarray*}
 \int_K  \vert \nabla _g p^g_{t/2}(y)\vert_{g}^{2}  d\mu_g(y)&= & \int -\Delta_gp^g_{t/2}(y) p^g_{t/2}(y)d\mu_g(y)\\
 &\leqslant & \| \Delta_gP^g_{t/4}p^g_{t/4}\|_2\|p^g_{t/2}\|_2 \leqslant \frac{4}{t} \|p^g_{t/4}\|_2 \|p^g_{t/2}\|_2
 \end{eqnarray*}
Recall that by \eqref{e.1.1} if $K$ is uniformly doubling with constant at most $D$, there is a constant $C(D)$ such that
 \begin{equation}\label{e.8.12}
 p_t^{g}(e)\leqslant  \frac{C \left( D \right)}{V_g(\sqrt{t})}.
 \end{equation}
 Moreover, by \eqref{e.DoublingComparison}, for any $0<a<1$
 \[
 p_{at}^{g}(e)\leqslant  \frac{C \left( D \right)}{D a^{\delta/2}V_g(\sqrt{t})}
 \]
 This yields
\[
\vert \nabla_g p_t^g(x)\vert_{g}\leqslant \frac{2}{\sqrt{t}}\left( p_{t}^{g}\left( e \right)\right)^{3/4}\left( p_{t/2}^{g}\left( e \right)\right)^{1/4}\leqslant \frac{C_1(D)}{t^{1/2} V_g(\sqrt{t})}.
\]
From \cite{Saloff-Coste2010b} or, more directly,  \cite[Theorem 4.11]{CoulhonSikora2008} (see also \cite{AuscherCoulhonDuongHofmann2004}) Equation \ref{e.2.5} follows.

Inequality \eqref{e.2.6} follows by integration. See, e.g., \cite[Lemma 5.2.13]{Saloff-CosteBook2002}.
\end{proof}

\begin{rem}\label{r.OpenQuestion}
We do not know if it is possible to prove statements \eqref{e.2.3} and \eqref{LiYau} with a uniform constant ($H_1$ or $H'_1$) over all $g \in \mathfrak{L}(K)$ in the case when $M=K$ is a compact Lie group, solely from the validity of Conjecture \ref{MC1}. In particular, we do not know if these statements hold uniformly for all left-invariant metrics on $\SU(2)$. There seems to be no reasons why they should not hold but the known techniques to attack these problems usually involve curvature.

In this direction we note that  the heat kernel lower bound in Theorem \ref{t.HeatKernelLower}  and \eqref{e.2.5} imply that there exist $C>0$ and  $b >1$ such that for all $x, t$
\begin{equation}\label{e.2.7}
|\nabla_g p_t(x)|_g \leqslant  Ct^{-1/2} p_{bt}(x).
\end{equation}
This is \eqref{e.2.3} for the heat kernel $p_{t}^{g} \left( x \right)$. Note that given \eqref{e.1.1} Equation \eqref{e.2.7} is equivalent to
\[
|\nabla_g p_t(x)| \leqslant   \frac{C}{\sqrt{t}V_g(\sqrt{t})}\exp(- b'
|x|_g^2/t).
\]
All the constants depend only on $D$ as follows from the proofs in  \cite{AuscherCoulhonDuongHofmann2004}.
\end{rem}
Finally, Theorem \ref{t.GradientEstimates} by  \cite{AuscherCoulhonDuongHofmann2004}  gives the following corollary regarding the Riesz transforms.

\begin{cor}[Uniform boundness of Riesz transforms] Assume $K$ is uniformly doubling with constant at most $D$. Then for all $1<p<\infty$  there are $c_{p}\left( D \right), C_{p}\left( D \right)$ such that
\[
c_{p}\left( D \right)\|\Delta_g^{-1/2}f\|_p \leqslant \||\nabla_g f|_g\|_p \leqslant
C_{p}\left( D \right)\|\Delta_g^{-1/2}f\|_p.
\]
\end{cor}

\subsection{Weyl counting function}\label{ss.WeylCount}
For a compact Riemannian manifold $\left( M, g \right)$, consider the Weyl spectral counting function
\[
\mathfrak W_{M, g}(s):=\#\{i :\lambda_{g, i} <s\},
\]
where $0=\lambda_{g, 0} < \lambda_{g} \leqslant \cdots\leqslant
\lambda_{g, i}\leqslant \cdots$ are the eigenvalues of $\Delta_g$ as
defined in \eqref{e.spectrum}.  The asymptotic behavior of this
function is described classically by Weyl's law (see \cite[p. 155]{ChavelBook1984}) as follows.
\begin{equation}
  \mathfrak W_{M, g}(s) \sim \frac{\omega_n}{(2\pi)^n} \mu_g(M) s^{n/2},
\end{equation}
where $\omega_n$ is the volume of the Euclidean $n$-ball.  However, even when $M=K$ is a compact connected Lie group, these asymptotics do not hold uniformly over all left-invariant metrics $g$ -- not even when $K = \mathbb{T}^n$ is a torus.

When $(M,g)$ is a compact homogeneous space, C.~Judge and R.~Lyons in \cite{JudgeLyons2017} have recently obtained the following uniform upper bound.
\begin{equation}
  \mathfrak{W}_{M,g}(s) \leqslant C \frac{\mu_g(M)}{V_g(s^{-1/2})},
\end{equation}
where $C$ is a universal constant.  If $M=K$ is a compact connected Lie group which is uniformly doubling, we
obtain a matching lower bound, uniformly over all left-invariant metrics.

\begin{pro}\label{weyl-upper}
Let $K$ be a connected compact real Lie group which is uniformly
doubling with constant at most $D$. Then there is a constant $c(D) > 0$,
depending only on $D$, such that for all $g \in \mathfrak{L}(K)$ we have
\begin{equation}\label{Weyl}
\mathfrak W_{K, g}(s) \geqslant c(D)\frac{\mu_{g}(K)}{V_{g}(s^{-1/2})}.
\end{equation}
\end{pro}

For a proof, see \cite[Th\'{e}or\`{e}me 7.1]{MaheuxSaloff-Coste1995};
an explicit statement is also given in \cite[Theorem 4.2]{SaloffCoste1995a}. The proof is based on the min-max
characterization of eigenvalues and a covering argument.  A matching
upper bound is also proved in \cite{MaheuxSaloff-Coste1995}, but with
a constant depending on $D$. Such a bound can be obtained another way
using the trace of the heat kernel, via \eqref{e.HeatKernelEigenvalue} and the heat kernel estimates of
Section \ref{ss.HeatKernelEst}. In \cite[Theorem 2]{FeffermanPhong1983} a similar statement is proved for individual sub-elliptic operators in $\mathbb{R}^{n}$, but without explicit
control of the constants involved in terms of doubling.

The bound \eqref{Weyl} is informative for $t \geqslant c
\operatorname{diam}_g^{-2}$ as discussed in Section
\ref{ss.spectral-gap}. Indeed, the spectral gap estimate in Theorem
\ref{spectral-gap} implies that  the step function $\mathfrak
W_{K,g}(t)$ equals $1$ on $[0, c \operatorname{diam}_{g}^{-2})$ for some $c=c(D)$, uniformly over all left-invariant metrics in $\mathfrak{L}(K)$.

When $K = \SU(2)$, Theorem \ref{thV} and Proposition \ref{diameter-bounds} yield detailed explicit estimates for $W_{\SU(2), g}(s)$ as follows.

\begin{cor}[Weyl counting function for $\SU(2)$] There are constants
  $0< C_{0}(D) \leqslant C_{\infty} < \infty$, with $C_0(D)$ depending
  only on $D$ and $C_\infty$ universal, such that for all $g = \in
  \mathcal{L}\left( \SU(2) \right)$ we have
\begin{equation}
C_{0}(D) f_{a_{1}, a_{2}, a_{3}}\left(t\right)  \leqslant
\mathfrak W_{\SU(2), g}(t) \leqslant
C_{\infty}f_{a_{1}, a_{2}, a_{3}}\left(t\right),
\end{equation}
where
\begin{equation*}
f_{a_{1}, a_{2}, a_{3}}\left(t\right)=
\begin{cases}
1 & \text{ if }  0 < t <  1/a_{2}^{2},
\\
a_2^{2} t & \text{ if }   1/a_{2}^{2} \leqslant t < 1/a_{1}^{2}
\\
a_1^2 a_2^{2} t^{2} & \text{ if }   1/a_{1}^{2} \leqslant t < a_3^{2}
/ a_{1}^{2}a_{2}^{2}
\\
a_1 a_{2} a_3 t^{3/2} & \text{ if } a_3^{2} / a_{1}^{2}a_{2}^{2} \leqslant t < \infty.
\end{cases}
\end{equation*}
Here $a_1, a_2, a_3$ are the parameters of $g$ as in Notation \ref{n.MetricsParameters}.
\end{cor}

\subsection{Heat kernel estimates: ergodicity}

Let $\mathbf V_g$ be the total Riemannian volume of  the given group $K$ under a Riemannian metric $g \in \mathfrak{L}(K)$, that is, $\mathbf V_g=\mu_{g}(K)$. It is well-known that the heat semigroup associated  to any given $g\in \mathfrak{L}(K)$ is ergodic and that $p^g_t(x)\longrightarrow  \mathbf V_g^{-1}$ as $t$ tends to infinity. As before let $\lambda_{g}$  by the lowest non-zero eigenvalue of the Laplacian $\Delta_g$ on $K$. We would like to describe this convergence to equilibrium in terms of the eigenvalue $\lambda_{g}$  in the case when $K$ is a uniformly doubling compact Lie group. For relevant results we refer to  \cite{Saloff-Coste1994a, Saloff-Coste2010b}. In what follows we set $\Vert f \Vert_{1}=\Vert f \Vert_{L^{1}\left( K, \mu_{g} \right)}$ and $\Vert f \Vert_{2}=\Vert f \Vert_{L^{2}\left( K, \mu_{g} \right)}$.

\begin{theo}
Let $K$ be a compact Lie group which is uniformly doubling with constant at most $D$.  For any $\epsilon>0$ there is a constant $C_\varepsilon(K)\in (0,\infty)$ such that for any metric $g\in \mathfrak{L}(K)$ we have
\[
\mathbf V_g \|p^g_t- \mathbf V_g^{-1}\|_1 \geqslant e^{-t\lambda_{g}} \hskip0.1in \text{ for all } t>0,
\]
and
\[
\mathbf V_g \|p^g_{t}-\mathbf V_g^{-1}\|_2^2 \leqslant  C_\varepsilon(K) e^{-2s\lambda_{g}}  \text{ for }  t \geqslant \varepsilon \operatorname{diam}_g^2 +s.
\]
Moreover, there are constants $c_i\in (0,\infty)$, $1\leqslant i \leqslant 4$, such that   for any metric $g\in \mathfrak{L}(K)$ we have
\[
\frac{c_1}{V_g(\sqrt{t})}e^{-2\lambda_{g} t} \leqslant
  \|p^g_{t}-\mathbf V_g^{-1}\|_2^2
 \leqslant  \frac{c_2}{V_g(\sqrt{t})}  e^{-2\lambda_{g} t}  \hskip0.1in \text{ for all } t>0.
 \]
 \end{theo}
\begin{rem} As we described in Section \ref{ss.spectral-gap}, under the hypothesis of this Theorem,  $\lambda_{g}$ is of order $\operatorname{diam}_g^{-2}$, uniformly over $\mathfrak{L}(K)$. Note also that, by definition,
$\mathbf V_g=V_g(\operatorname{diam}_{g}/2)$.   Further, for any function $f\in L^2(K, \mu_{g})$, $\| f \|_1^{2}\leqslant  \mathbf V_g \|f\|_2^{2}$.
\end{rem}
\begin{proof}
Let $\varphi$ be an eigenfunction of $\Delta_g$ associated with the lowest non-zero eigenvalue $\lambda_{1, g}$ and normalized by $\max_x |\varphi(x)| =\varphi(e)=1$ (such a normalization is always possible by translation in $K$ and multiplication by a constant). Then the lower $L^1$-bound follows from
\[
\|p^g_t-\mathbf V_g^{-1}\|_1\geqslant \int \left(p^g_t(x)-\mathbf V_g^{-1}\right)\varphi(x)d\mu_g(x)=e^{-\lambda_{g}t},
\]
where we used the fact that $\varphi \perp 1$ and
\[
\int_{K} p^g_t(x)\varphi(x)d\mu_g(x)=\left(P_{t}\varphi\right)\left( e \right)=e^{-\lambda_{g}t}\varphi(e).
\]

For the two-sided $L^2$-estimate, first observe that for any constant $C>$ by \eqref{e.8.1} we have
\[
\int_{K} \vert p^g_t(x) -C \vert^{2} d\mu_{g}\left(x\right)=p^g_{2t}(e)-2C+C^{2}\mathbf V_g,
\]
and so by \eqref{e.HeatKernelEigenvalue}
\[
\mathbf V_g\|p^g_{t}-\mathbf V_g^{-1}\|_2^2= \mathbf V_g p^g_{2t}(e) -1 =\sum_{i=1}^\infty e^{-2t \lambda_{g,i}},
\]
where  $\lambda_{g,i}$  are eigenvalues of $\Delta_g$ as defined in
\eqref{e.spectrum}. For the lower bound, noting that $e^{-2t
  \lambda_g} \leqslant 1$ we have
\begin{align*}
\mathbf V_g p^g_{2t}(e) -1 =
e^{-2t \lambda_g} + \sum_{i=2}^\infty e^{-2t \lambda_{g,i}}
&\geqslant \frac{1}{2}e^{-2t \lambda_g}\left(1 + e^{-2t \lambda_g}\right) +
\sum_{i=2}^\infty e^{-2t \lambda_{g,i}} \\
&\geqslant \frac{1}{2} e^{-2t\lambda_{g}} \left(1 + e^{-2t
  \lambda_g}\right) + \frac{1}{2}e^{-2t \lambda_g} \sum_{i=2}^\infty
e^{-2t \lambda_{g,i}} \\
&= \frac{1}{2} e^{-2t \lambda_g} \mathbf V_g p^g_{2t}(e).
\end{align*}
By  Theorem \ref{t.HeatKernelLower}, this gives the desired lower bound.

For the upper bound, write
\begin{align*}
\mathbf V_g p^g_{2t}(e) -1 = \sum_{i=1}^\infty e^{-2t
  \lambda_{g,i}} &=\sum_{\lambda_{g,i} \le 2 \lambda_g} e^{-2t
  \lambda_{g,i}}+\sum_{\lambda_{g,i} > 2 \lambda_g} e^{-2t \lambda_{g,i}}
\\
&\leqslant  e^{-2t\lambda_{g}} \left( \mathfrak W_g(2\lambda_{g}) + \sum_{\lambda_{g,i} > 2 \lambda_g} e^{-t\lambda_{g,i}}\right)
\\
&\leqslant e^{-2t\lambda_g} \left(  \mathfrak W_g(2\lambda_{g}) +  \mathbf V_g p^g_{t}(e)\right).
\end{align*}
The upper bound in \eqref{Weyl} and the upper bound $\lambda_g \leqslant C_1(D) \operatorname{diam}_g^{-2}$ in Theorem \ref{spectral-gap} yield
$\mathfrak W_g(2\lambda_g)\leqslant C_2(D)$. This, together
with \eqref{e.HeatKernelUpper}, gives
\[
\mathbf V_g p^g_{2t}(e) -1 \leqslant \frac{C_3(D) \mathbf V_g}{V_g(\sqrt{t})} e^{-2t\lambda_g}
\]
as desired.
\end{proof}

\subsection{Infinite products}\label{s.InfiniteProducts}
Let $\left\{ K_{i} \right\}_{i=1}^{\infty}$ be a sequence of compact connected Lie groups, each equipped with the Haar probability measure $\mu_{i}$. Consider the compact group
 \[
 K^\infty:=\prod_{i=1}^{\infty} K_i.
  \]

Note that this includes the case when $K_{i}=K$ for all $i$. Suppose each $K_i$ is equipped with a Riemannian metric $g_i \in \mathfrak{L}(K_{i})$; from now on by $K_i$ we denote $(K_{i}, g_i)$,
and by $\mathbf g$ we denote the sequence of metrics $\left\{g_i\right\}_{i=1}^{\infty}$. Note that the Riemannian volume measure $\mu_{g_i}$ is just a rescaling of $\mu_i$, so $D(K_i, d_{g_i},
\mu_{g_i}) = D(K_i, d_{g_i}, \mu_i)$.  We endow $K^\infty$ with its Haar probability measure $\mu$ which is the product of the Haar measures $\mu_{i}$. For background on this setting see  \cite{HeyerBook1977, BendikovSaloffCoste2001a}.

The space of \emph{cylinder functions}, i.e. smooth functions depending on only finitely many coordinates, is dense in $L^2(K^\infty,\mu)$.  For a cylinder function $f$, set
 \[
 \mathcal E_{\mathbf g}(f,f):=\int _{K^\infty} \sum_{i=1}  g_i(\nabla_{ g_i} f,\nabla_{g_i}f)d\mu.
 \]
 The quadratic form $ \mathcal E_{\mathbf g}$ is closable and its closure is a strictly local regular Dirichlet form associated to a self-adjoint Markov semigroup $H^{\mathbf g}_t$. It is a convolution semigroup on $K^\infty$ associated with a convolution semigroup of symmetric measures $\nu^{\mathbf g}_t$, i.e.
 \[
 H^{\mathbf g}_tf (x)=\int f(xy)d\nu^{\mathbf g}_t(y), \hskip0.1in t>0.
 \]
 For each metric $g_i$, we let $\gamma_i:=\lambda_{1,i}$ be the second smallest eigenvalue of the operator $-\Delta_{i}$, where $\Delta_{i}$ is the Laplace--Beltrami operator on $K_i$.

Denote by $t_{A}$ to be the infimum of all times $t > 0$ at which the measure $\nu^{\mathbf g}_t$ is
absolutely continuous with respect to the Haar measure $\mu$. Note that if this property holds at time $t$, it also holds at all later times.

The following are special cases of more general open problems considered in  \cite[Section 2]{Saloff-Coste2010b}. Is is true that $\nu^{\mathbf g}_t$ is singular with respect to the Haar measure $\mu$ for all time $t < t_{A}$? Is it true that for all $t > t_{A}$, the density $f^{\mathbf g}_t$ of the measure $\nu^{\mathbf g}_t$ with respect to $\mu$ is in $L^{2}\left( K^{\infty}, \mu \right)$? Is
it true that if $t_{A}=0$, then $f^{\mathbf g}_t$  has a continuous representative?

\begin{pro} Assume there exists a constant $D$ such that for any $i=1, 2, ....$ and any $g_{i} \in \mathfrak{L}\left( K_{i} \right)$ we have $D\left( K_{i}, d_{g_{i}}, \mu_{i}\right) \leqslant D$. Denote
\[
t_* :=\inf\{t: \sum_1^\infty e^{-2t\gamma_i} <\infty\}.
\]
Then the following properties hold:

\begin{itemize}
\item The measure $\nu^{\mathbf g}_t$ is absolutely continuous with respect to the Haar measure $\mu$ for $t>t_*$ whereas $\nu^{\mathbf g}_t$  has no absolutely continuous part with respect to $\mu$ for $0<t<t_*$;

\item Furthermore, for all $t>t_*$, the density  $\frac{d\nu^{\mathbf g}_t}{d\mu}$ is
in $L^2(K^{\infty},\mu)$. It is unbounded for $t_*<t<2t_*$, and it is bounded and continuous for $t>2t_*$;

\item In particular,  if $t_*=0$, the semigroup $H^{\mathbf g}_t$
  admits a continuous convolution kernel for all times $t>0$.
  \end{itemize}
\end{pro}

\begin{proof}
This follows from \eqref{e.1.0}, \eqref{e.1.1} Proposition \ref{p.2.3} and \cite[Theorems 3.1, 4.1, 4.2]{BendikovSaloff-Coste1997a}.
\end{proof}
In particular, by Theorem \ref{thM} and the similar result for tori, these properties hold when $K_{i} \in \left\{ \SU(2), \mathbb{T}, \mathbb{T}^{2}, \dots, \mathbb{T}^{n}\right\}$.

\section{Connections to sub-Riemannian geometry}\label{sub-riemannian}

\newcommand{\Lsub}{\mathfrak{L}_{\mbox{\tiny sub}}}

We have focused this paper on Riemannian geometry, but in fact our
results carry over to sub-Riemannian geometry as well.  In this
section, we make those connections explicit.  We briefly review the
relevant definitions as they apply to Lie groups; we refer to
\cite{MontgomeryBook2002} for a discussion of sub-Riemannian geometry
in a general context.

On a connected Lie group $K$, a left-invariant sub-Riemannian geometry
is determined by a choice of a linear subspace $H \subset
\mathfrak{k}$ of the Lie algebra, and a Euclidean inner product $g$ on
$H$.  Let $\Lsub(K)$ denote the set of all such pairs $(H,g)$; by
abuse of notation, we will refer to such a pair simply by $g$.  It is
also common to view $g$ as an extended quadratic form on
$\mathfrak{k}$, where $g(v,w) = \infty$ unless $v,w \in H$.

By left translation, $H$ extends to a left-invariant distribution
$\mathcal{H} \subset TK$ with $\mathcal{H}_e = H$, and $g$ extends to
a left-invariant sub-Riemannian metric, still called $g$, on
$\mathcal{H}$ (or an extended quadratic form on $TK$).

The geometry $(H,g)$ satisfies the H\"ormander bracket generating
condition iff $H$ generates the Lie algebra $\mathfrak{k}$; let
$\Lsub^*(K) \subset \Lsub(K)$ denote the set of such geometries.  Note
that for $K = \SU(2)$, this happens iff $\dim H \geqslant 2$, since
the Lie algebra $\su(2)$ is generated by any two linearly independent
elements.  When $H = \mathfrak{k}$ we recover the left-invariant
Riemannian geometries $\mathfrak{L}(K)$.

To any $g\in \Lsub(K)$ is associated a length structure giving finite
length to continuous piecewise smooth curves that stay tangent to
$\mathcal{H}$ (these are called horizontal curves).  The
left-invariant Carnot--Carath\'eorody (pseudo)-distance $d_g(x,y)$ is
defined as the infimum of the lengths of horizontal curves joining $x$
to $y$ in $K$, where $d_g(x,y)=\infty$ if no such curve
exists.  By the Chow--Rashevskii theorem \cite[Theorems
  2.1.2 and 2.1.3]{MontgomeryBook2002}, if $g\in \Lsub^*(K)$ then $d_g(x,y)$
is finite for any pair $x,y\in K$, so that $d_g$ is a genuine
distance, and moreover the topology induced by $d_g$ coincides with
the manifold topology of $K$.

Each sub-Riemannian geometry $(H,g) \in \Lsub(K)$ is also associated with
a canonical left-invariant sub-Laplacian $\Delta_g$, which may be defined by
\begin{equation}
  \Delta_g = - \sum_{i=1}^{k} \widetilde{u_i}^2
\end{equation}
where $k = \dim H$, $\{u_i : 1 \leqslant i \leqslant k \}$ is a
$g$-orthonormal basis for $H$, and $\{\widetilde{u_i}\}$ are the
corresponding left-invariant vector fields.  This definition is
independent of the basis chosen.  The operator $\Delta_g$ is
hypoelliptic iff $g \in \Lsub^*(K)$, and when $g$ is Riemannian ($H =
\mathfrak{k}$) we recover the Laplace--Beltrami operator.

Likewise, for $f \in C^\infty(K)$, we have the left-invariant
sub-gradient $\nabla_g f$ which is a smooth section of $\mathcal{H}$
defined by
\begin{equation}
  \nabla_g f = \sum_{i=1}^k (\widetilde{u_i} f) \widetilde{u_i}.
\end{equation}
In particular, we have
\begin{equation}
  |\nabla_g f|^2 := g(\nabla_g f, \nabla_g f) = \sum_{i=1}^k
  |\widetilde{u_i} f|^2.
\end{equation}
When $g$ is Riemannian this is the usual Riemannian gradient.

In the case $K = \SU(2)$, a sub-Riemannian metric $g \in
\Lsub(\SU(2))$ can be diagonalized by a standard Milnor basis, in the
same way as in Lemma \ref{milnor-diagonalize} for Riemannian metrics.

\begin{pro}\label{milnor-diagonalize-sub}
  Let $(H,g)\in \mathcal \Lsub(\SU(2))$, with $\dim H = k$. There
  exists a standard Milnor basis $\{e_1,e_2,e_3\}$ and an ordered
  triplet of extended non-negative reals $0 <  a_1\leqslant a_2\leqslant a_3\leqslant
  \infty$ such that $H = \operatorname{span}\{ e_i : 1 \leqslant i \leqslant k\}$ and $g(e_i,
  e_j) = a_i^2 \delta_{ij}$ for $1 \leqslant i,j \leqslant k$.  We take $a_i =
  \infty$ for $i > k$.
\end{pro}

\begin{proof}
  The case $k=0$ is trivial (any standard Milnor basis will do), and
  $k=3$ is Lemma \ref{milnor-diagonalize}.

  For $k=2$, let $\{v_1, v_2\}$ be a $g$-orthonormal basis for $H$,
  and set $v_3 = [v_1, v_2]$.  Observe that $v_3 \notin H$; indeed,
  under the invariant inner product given by the negative Killing
  form, $v_3$ is orthogonal to both $v_1, v_2$.  Let $g'$ be the
  Euclidean inner product on $\su(2)$ which makes $v_1, v_2, v_3$
  orthonormal, and define $\times, L$ with respect to $g'$ as in the
  proof of Lemma \ref{milnor-diagonalize}, choosing $\times$ so that
  $v_1 \times v_2 = v_3$.  Note that $v_3$ is an
  eigenvector of $L$ (with eigenvalue $1$), since $L(v_3) = L(v_1
  \times v_2) = [v_1, v_2] = v_3$.  So if $\{w_1, w_2, w_3\}$ is a
  $g'$-orthonormal basis of eigenvectors for $L$, where we let $w_3 =
  v_3$, then necessarily $w_1, w_2 \in H$ and they are $g$-orthonormal.
  Proceeding as in Lemma \ref{milnor-diagonalize}, there is a standard
  Milnor basis $\{e_1, e_2, e_3\}$ where $e_i$ is a scalar multiple of
  $w_i$, and in particular $e_1, e_2 \in H$ and they are
  $g$-orthogonal.

  For $k=1$, let $v_1$ span $H$, choose $v_2 \notin H$ arbitrarily,
  and proceed as in the previous case.  We obtain a standard Milnor
  basis $\{e_1, e_2, e_3\}$ where $\operatorname{span}\{e_1, e_2\} =
  \operatorname{span}\{v_1, v_2\}$.  In particular there is some
  $\theta \in \mathbb{R}$ such that $v_1$ is a scalar multiple of
  $\cos(\theta) e_1 + \sin(\theta) e_2$, and then
  \begin{equation*}
    \{\cos(\theta) e_1 + \sin(\theta) e_2, \sin(\theta) e_1 -
    \cos(\theta)e_2, e_3\}
  \end{equation*}
  is the desired standard Milnor basis, as in Example \ref{milnor-transformations}.
\end{proof}

Thus, as in Corollary \ref{isometry}, the left-invariant
sub-Riemannian geometries $g \in \Lsub(\SU(2))$ are given, up to
isometry, by the geometries $g_{(a_1, a_2, a_3)}$, where the $a_i$ are
allowed to take the value $\infty$.  In fact, these geometries arise
as the limits of the Riemannian geometries $g_{(a_1, a_2, a_3)}$ where
the $a_i$ are finite.  The ``standard'' sub-Riemannian metric
commonly encountered in the literature
(e.g. \cite{BaudoinBonnefont2009,BaudoinGarofalo2017}) corresponds to
$g_{(1, 1,\infty)}$, but we stress that this is just one element of the
infinite family $\Lsub^*(\SU(2))$.

\begin{lem}\label{lem-conv1}
  Given $g = g_{(a_1, a_2, a_3)} \in \Lsub(\SU(2))$, where $0 < a_1
  \leqslant a_2 \leqslant a_3 \leqslant \infty$, and $\epsilon > 0$,
  let $a_{\epsilon,i} = \min(a_i, \epsilon^{-1})$, and set $g_{\epsilon} =
  g_{(a_{\epsilon,1}, a_{\epsilon,2}, a_{\epsilon, 3})} \in
  \mathfrak{L}(\SU(2))$.  Then for any $x,y \in \SU(2)$ we have
  $d_g(x,y) = \lim_{\epsilon \to 0} d_{g_\epsilon}(x,y)$.
\end{lem}

\begin{proof}
By left invariance, it suffices to consider $d_g(e,x)$ where $x \ne e$.

If $d_g(e,x) < \infty$, the result follows by the argument in \cite[Proposition 3.1]{JerisonSanchezCalle1987} for the distance $\alpha_L$.  In particular, this covers all cases when $a_2 < \infty$ (so that $\dim H$ is $2$ or $3$).

  In the trivial case of $g_{(\infty, \infty, \infty)}$, where $\dim H
  = 0$, we have $d_g(e,x) = \infty$ for all $x \ne e$, and we simply
  note that $d_{g_\epsilon}(e,x) = d_{(\epsilon^{-1}, \epsilon^{-1},
    \epsilon^{-1})}(e,x) = \epsilon^{-1} d_{(1,1,1)}(e,x) \to \infty$
  as $\epsilon \to 0$.

  The remaining case is where $g = g_{(a_1, \infty, \infty)}$, with
  $a_1 < \infty$ (so that $\dim H = 1$) and $d_g(e,x) = \infty$.  Let
  $S = \{\exp(s \hat{e}_1) : s \in \mathbb{R}\}$ be the circle
  subgroup defined in the proof of Proposition \ref{diameter-bounds}.
  If $x \in S$, so that $x = \exp(T \hat{e}_1)$ for some $T$, then
  $\gamma(t) = \exp(t \hat{e}_1)$, $0 \leqslant t \leqslant T$ is a
  finite-length horizontal curve joining $e$ to $x$, and thus
  $d_g(e,x) < \infty$.  So suppose $x \notin S$.  As shown in the
  proof of Proposition \ref{diameter-bounds}, we have
  $d_{(0,1,1)}(e,x) > 0$.  Hence for all $\epsilon \leqslant
  \min(a_1^{-1}, 1)$ we have
  \begin{equation*}
    0 < d_{(0,1,1)}(e,x) \leqslant d_{(\epsilon a_1, 1, 1)}(e,x) =
    \epsilon d_{g_{\epsilon}}(e,x)
  \end{equation*}
  which implies that $d_{g_{\epsilon}}(e,x) \to \infty$.
\end{proof}

\begin{cor}\label{subriem-doubling}
  The family of metric measure spaces
  $$\{ (\SU(2), d_g, \mu_0) : g \in
  \Lsub^*(\SU(2))\}$$
  is uniformly volume doubling with the same
  constant $D$ as in Theorem \ref{thM}.
\end{cor}

\begin{proof}
  By the previous lemma, the closed ball $\bar{B}_g(r)$ equals the
  decreasing intersection $\bigcap_n B_{g_{1/n}}(r)$.  The
  sub-Riemannian spheres have measure zero \cite[Proposition
    4.3]{RiffordTrelat2005}, so we have $\mu_0(B_g(r)) =
  \mu_0(\bar{B}_g(r)) = \lim_{\epsilon \to 0}
  \mu_0(B_{g_\epsilon}(r))$, and by Theorem \ref{thM} each
  $g_\epsilon$ is volume doubling with constant at most $D$, so the
  result follows.
\end{proof}

\begin{cor}
  For all $g_{(a_1, a_2, a_3)} \in \Lsub^*(\SU(2))$, where we allow
  $a_3 = \infty$, the volume $V_{(a_1, a_2, a_3)}(r)$ is comparable to
  $\overline{V}_{(a_1, a_2, a_3)}(r)$ as defined in
  \eqref{Vbar-a-def}, uniformly in $a_1, a_2, a_3, r$.
\end{cor}

Note that for $a_3 = \infty$, the ``Euclidean'' regime, where volume
scales as $r^3$, becomes empty, and for very small $r$, the volume
scales as $r^4$ instead.  This matches the Heisenberg behavior and
corresponds to the fact that such a sub-Riemannian geometry has
Hausdorff dimension 4.

\begin{rem}
  The preceding corollaries may also be proved directly, instead of by
  approximating sub-Riemannian geometries by Riemannian geometries.
  Indeed, the proofs in Sections
  \ref{s.EuclideanRegime}--\ref{s.Combining} go through without change
  if $a_3 = \infty$.  (Note that Section \ref{s.EuclideanRegime}, the
  Euclidean regime, becomes vacuous in that case.)
\end{rem}

The results in Section \ref{s.Consequences} concerning the spectral
gap $\lambda_g$, the heat kernel $p_t^g$, the eigenvalues
$\lambda_{g,i}$ and the Weyl counting function $\mathfrak{W}_g$ all
extend uniformly to sub-Riemannian geometries $g \in \Lsub^*(\SU(2))$,
with $\Delta_g$, $|\nabla_g f|^2$ redefined as above.  It is only
necessary to adjust the statements to replace all instances of
$\mu_g$ by $\mu_0$, since sub-Riemannian geometries do not admit a
Riemannian volume, and scale appropriately.  In particular, in this
context the heat kernel $p_t^g$ should be viewed as an integral kernel
with respect to $\mu_0$.  The proofs need not be carried out by
passing to the limit in the Riemannian statement; instead, the results
follow because they are general consequences of uniform doubling and
the uniform Poincar\'e inequality, for which the proof cited in
Section \ref{a.Poincare} goes through without change in the
sub-Riemannian setting.

\begin{rem}
  One may also study the degenerate sub-Riemannian geometries, though
  this is more complicated because their topologies are not well
  behaved.  For instance, with $g = g_{(a_1, \infty, \infty)}$, the
  $\infty$-metric space $(\SU(2), d_g)$ has uncountably many connected
  components, which are the left cosets of the one-dimensional
  subgroup $S = \{\exp(s e_1) : s \in \mathbb{R}\}$, all isometric to
  $S^1$ and at pairwise distance infinity from one another.  In
  particular, every ball of this metric has Haar measure zero, so
  statements about volume growth are not sensible.  However, if we fix
  a sufficiently large $R$, then for all small $\epsilon$ the ball
  $B_{g_\epsilon}(R)$ is comparable to $B_{g_{\epsilon}}(S,R)$; by
  arguments similar to Proposition \ref{collapse-upper}, one may see
  that $\mu_0(B_{g_\epsilon}(R)) \approx \epsilon^2 R^2$.  On the
  other hand, \eqref{Vbar-a-def} gives
  \begin{equation*}
  \mu_0(B_{g_{\epsilon}}(r)) \approx \begin{cases} a_1^{-1}
      \epsilon^2 r^3, & 0 \leqslant r \leqslant a_1 \\
      \epsilon^2 r^2, & a_1 \leqslant r \leqslant \epsilon^{-1} \\
      1, & r \geqslant \epsilon^{-1}.
      \end{cases}
  \end{equation*}
  As $\epsilon \to 0$, the ball $B_{g_{\epsilon}}(R)$ collapses to
  $S$, and we have
  \begin{equation}\label{circle-vol}
    \frac{\mu_0(B_{g_{\epsilon}}(r))}{\mu_0(B_{g_{\epsilon}}(R))}
    \approx
    \begin{cases}
      \frac{r}{a_1}, & r \leqslant a_1 \\
      1, & r \geqslant a_1.
    \end{cases}
  \end{equation}
  If we consider the circle $S^1$ as a Lie group equipped with its own
  normalized Haar measure $\mu_{S^1}$ and the metric $g = g_{a_1}$
  which is the $a_1$-scaling of the unique left-invariant Riemannian
  metric on $S^1$, we can observe that \eqref{circle-vol} is
  comparable to the volume $\mu_{S^1}(B_g(r))$ of a ball in $S^1$.  In
  particular, we recover the (trivial) fact that left-invariant
  Riemannian geometries on $S^1$ are uniformly volume doubling.  This
  is perhaps not so interesting in our present context, but the idea
  of considering degenerate sub-Riemannian geometries may yield more
  useful insights when replacing $\SU(2)$ with other compact connected
  Lie groups $K$.
\end{rem}

\begin{acknowledgement}
The third author would like to thank Dominique Bakry for inspiring discussions over many years.  The authors are grateful for helpful and motivating conversations with Iddo Ben-Ari, Bruce K.~Driver, Nicolas Juillet, Emilio A.~Lauret, and K.-T.~Sturm.
\end{acknowledgement}

\bibliographystyle{amsplain}
\def\cprime{$'$}
\providecommand{\bysame}{\leavevmode\hbox to3em{\hrulefill}\thinspace}
\providecommand{\MR}{\relax\ifhmode\unskip\space\fi MR }
\providecommand{\MRhref}[2]{%
  \href{http://www.ams.org/mathscinet-getitem?mr=#1}{#2}
}
\providecommand{\href}[2]{#2}

\end{document}